\documentclass[11pt]{amsart}
\usepackage{amsmath,amssymb,amsfonts}
\usepackage{enumerate}
\usepackage{amscd, amssymb, latexsym, amsmath, amscd}
\usepackage[all]{xy}
\usepackage{pb-diagram}
\usepackage{picture}
\usepackage{multicol,lipsum}

\newtheorem{theorem}{Theorem}[section]
\newtheorem{thm}[equation]{Theorem}
\newtheorem{prop}[equation]{Proposition}
\newtheorem{cor}[equation]{Corollary}
\newtheorem{con}[equation]{Conjecture}

\numberwithin{equation}{section}

\newcommand{\Q}{\mathbb Q}
\newcommand{\Z}{\mathbb Z}
\newcommand{\R}{\mathbb R}
\newcommand{\C}{\mathbb C}

\newcommand{\A}{\mathbb A}

\def\Hom{{\rm Hom}}

\def\SL{{\rm SL}}
\def\disc{{\rm disc}}

\def\PGSp{{\rm PGSp}}

\def\Sp{{\rm Sp}}

\def\SU{{\rm SU}}
\def\PU{{\rm PU}}
\def\U{{\rm U}}

\def\GL{{\rm GL}}
\def\PGL{{\rm PGL}}

\def\SO{{\rm SO}}

\def\Sp{{\rm Sp}}
\def\Mp{{\rm Mp}}

\def\O{{\rm O}}

\usepackage[all]{xy}

\def\A{{\mathbb A}}

\def\R{{\mathbb R}}

\def\Z{{\mathbb Z}}

\def\C{{\bf C}}

\def\C{{\mathbb C}}

\begin{document}

\title{Automorphic Forms and the Theta Correspondence} 
 \author{Wee Teck Gan}
 

\maketitle

\section{\bf Lecture 1: The Ramanujan Conjecture}

In the first lecture, we shall recall the Ramanujan conjecture for classical modular forms and its reformulation in the language of cuspidal automorphic representations of $\PGL_2$. 
For more details on this transition and reformulation, please take a look at \cite[\S 4.27 and \S 4.28]{G1}.
This reformulation allows one to readily generalize the conjecture to the setting of cuspidal automorphic representations of general connected reductive group $G$ over a number field $k$.   We will then discuss the unitary group analog of a construction of Roger Howe and Piatetski-Shapiro \cite{HPS} which gives a definitive counterexample to the extended Ramanujan conjecture for the unitary group $\U_3$ in three variables. Their original construction gave a counterexample on the group $\Sp_4$, but we will use the same idea to produce a counterexample on $\U_3$    via the method of theta correspondence. The unitary case we discussed here was  actually considered by Gelbart-Rogawski in \cite{GR}. 

\vskip 5pt

\subsection{\bf The Ramanujan conjecture}

About a century ago, Ramanujan considered the following power series of $q$
\[ \Delta(q) = q \cdot \prod_{n \geq 1} (1-q^n)^{24}. \]
Expanding this out formally, we have:
\[
\Delta(q) = \sum_{n >0} \tau(n) q^n = q- 24q^2+...\]
Ramanujan  conjectured that for all primes $p$,  
\[  |\tau(p)| \leq 2\cdot p^{11/2}. \]
 This is the Ramanujan conjecture in question.
  More generally, for any  holomorphic cuspidal Hecke eigenform $\phi$  of weight $k$ (and level $1$)  on the upper half plane $\mathfrak{h}$, with Fourier coefficients $\{a_n(\phi)\}_{n \geq 1}$,  the Ramanujan-Petersson conjecture asserts that
\[  |a_p(\phi)| \leq 2 p^{(k-1)/2}  \quad \text{for all primes $p$.} \]
For Hecke eigenforms, the Fourier coefficients $a_p(\phi)$ are also the eigenvalues of the Hecke operator $T_p$. Hence, the Ramanujan conjecture concerns bounds on cuspidal Hecke eigenvalues. It was eventually proved by Deligne \cite{D} as a consequence of his proof of the Weil's conjectures.
\vskip 5pt

\subsection{\bf Cuspidal automorphic representations}
The classical theory of modular forms can now be subsumed in a representation theoretic setting. The details of this transition can be found in \cite[\S 3 and \S 4]{G1}. Let us briefly recall this. 
\vskip 5pt

Let 
\vskip 5pt
\begin{itemize}
\item $k$ be a number field with ring of adeles $\A = \prod_v' k_v$, which is a restricted direct product of the local completions $k_v$ for all places $v$ of $k$.
\item $G$ be a connected reductive linear algebraic group over $k$; for simplicity we may take $G$ to be semisimple;
\item  $\rho: G \hookrightarrow \GL_n$ be a faithful algebraic representation over $k$, giving rise to a system of open compact subgroups $K_v = \rho^{-1}(\GL_n(\mathcal{O}_v)) \subset G(k_v)$ for almost all $v$, where $\mathcal{O}_v$ is the ring of integers of $k_v$; moreover, for almost all $v$, $K_v$ is hyperspecial.
\item$G(\A) = \prod_v' G(k_v)$ be the adelic group, which is a restricted direct product of $G(k_v)$ relative to the family $\{K_v \}$ of open compact subgroups for almost all $v$;
\item $[G] = G(k) \backslash G(\A)$ be the automorphic quotient; the locally compact group $G(\A)$ acts on $[G]$ by right translation and there is a $G(\A)$-invariant measure (unique up to scaling).
\end{itemize}
For the theory of classical modular forms, one is taking $k = \Q$ and $G = \PGL_2$. Using the natural identification
\[ \SL_2(\Z) \backslash  \mathfrak{h} \cong \PGL_2(\Z) \backslash  \PGL_2(\R) / \O_2(\R) \cong \PGL_2(\Q) \backslash \PGL_2(\A)/ \O_2(\R) \cdot  \prod_p \PGL_2(\Z_p),  \]
a classical modular form $\phi$ on $\mathfrak{h}$ corresponds to a function 
\[  f: \PGL_2(\Q) \backslash \PGL_2(\A) \longrightarrow \C \]
defined by
\[  f(g)  = (\phi|_k g) ( \sqrt{-1}). \]
Replacing $\phi$ by $f$ allows one to extend the notion of modular forms to the setting of general reductive groups $G$.

\vskip 5pt
More precisely, an automoprhic form on $G$ is a function 
\[  f: [G] \longrightarrow \C \]
satisfying some regularity and finiteness properties:
\begin{itemize}
\item $f$ is smooth
\item $f$ is right $K_f$-finite (where $K_f  = \prod_{v < \infty} K_v$)
\item $f$ is of uniform moderate growth
\item $f$ is $Z(\mathfrak{g})$-finite.
\end{itemize}
It is not important for us to know precisely the meaning of the above properties.  
\vskip 5pt

Let us denote the space of automorphic forms on $G$ by $\mathcal{A}(G)$.  
The group $G(\A)$ acts on the vector space $\mathcal{A}(G)$ by right translation.
An {\em automorphic representation} $\pi$ of $G$ is by definition an irreducible subquotient of the $G(\A)$-module $\mathcal{A}(G)$.
As an irreducible abstract representation of $G(\A) = \prod'_v G(k_v)$, $\pi$ is of the form:
\[  \pi = \otimes_v' \pi_v, \]
a restricted tensor product of irreducible smooth representations $\pi_v$ of $G(k_v)$. In particular, for almost all $v$, $\pi_v^{K_v} \ne 0$; we say that $\pi_v$ is $K_v$-spherical or $K_v$-unramified.   It is known that $\dim \pi_v^{K_v} =1$ if $\pi_v$ is $K_v$-unramified.

\vskip 5pt
\noindent{\bf Definition:}  An automorphic form $f$ on $G$ is called a {\bf cusp form} if, for any 
parabolic $k$-subgroup $P = MN$ of $G$, the $N$-constant term 
\[  f_N(g) = \int_{N(k) \backslash N(\A)} f(ng) \, dn \]
is zero as a function on $G(\A)$.
\vskip 5pt
Let
\[  \mathcal{A}_{cusp}(G) \subset \mathcal{A}(G) \]
be the subspace of cusp forms; it is a $G(\A)$-submodule (potentially 0). It turns out that, if $G$ is semisimple or has anisotropic center,  a cusp form $f$ is rapidly decreasing (modulo center)  as a function on $[G]$   and hence is square-integrable on $[G]$, i.e.
\[  \int_{[G]} |f(g)|^2 \, dg < \infty. \]
If we let $\mathcal{A}_2(G)$ denote the $G(\A)$-submodule of square-integrable automorphic forms (clearly nonzero), then 
\[  \mathcal{A}_{cusp}(G) \subset \mathcal{A}_2(G) \subset \mathcal{A}(G). \]
Moreover, each $G(\A)$-submodule $\mathcal{A}_{cusp}(G)$ and $\mathcal{A}_2(G)$ decomposes as a direct sum
\[  \mathcal{A}_{cusp}(G) = \bigoplus_{\pi} m_{cusp}(\pi) \cdot \pi \quad \text{and} \quad \mathcal{A}_2(G) = \bigoplus_{\pi} m_2(\pi)  \cdot \pi  \]
with finite multiplicities (this is not obvious). An irreducible summand $\pi$ of $\mathcal{A}_{cusp}(G)$ is called a cuspidal automorphic representation, whereas one of $\mathcal{A}_2(G)$ is called a square-integrable automorphic representation.

\vskip 10pt
\subsection{\bf Classification of unramified representations}
For  an abstract irreducible representation $\pi \cong \otimes'_v \pi_v$, we have mentioned that for almost all $v$, $\pi_v$ is $K_v$-unramified, with $K_v \subset G(k_v)$ a so-called hyperspecial maximal compact subgroup. Such $K_v$-unramified representations can be classified. We recall this classification briefly; the reader can take a look at 
\cite[\S 4.22-4.26]{G1}.
\vskip 5pt

 For almost all such $v$, $G$ is quasi-split over $k_v$ and hence possesses a Borel subgroup $B_v = T_v \cdot U_v$ defined over $k_v$. Hence, $T_v$ is a maximal torus over $k_v$.
 For a character $\chi_v: T_v = T_v(k_v) \rightarrow \C^{\times}$, one can form the (normalized) parabolically induced representation
 \[ I(\chi_v) =  {\rm Ind}_{B_v}^{G_v} \chi_v,  \]
 consisting of functions $\phi: G_v \rightarrow \C$ satisfying
 \[  f(utg) = \chi_v(t) \cdot \delta_{B_v}(t)^{1/2} \cdot f(g) \quad \text{for $u \in U_v$, $t \in T_v$ and $g \in G_v$,}  \]
  where $\delta_{B_v}$ is the modulus character defined by
\[  \delta_{B_v}(g) = | \det({\rm Ad}(g) | {\rm Lie}(U_v))| \]
and where the action of   $G_v$ is by right translation. The representation $I(\chi_v)$ is called a principal series representation.

\vskip 5pt

If  $\chi_v$  is an unramified character of $T(k_v)$, i.e., $\chi_v$ is trivial on $T(k_v) \cap K_v$, we call $I(\chi_v)$ an unramified principal series representation. The following proposition summarizes the classification of $K_v$-unramified representations:
\vskip 5pt

\begin{prop}
(i)  An  unramified $I(\chi_v)$ is of finite length and contains exactly one irreducible subquotient $\pi(\chi_v)$ which is $K_v$-unramified.  Moreover, 
$\pi(\chi_v) \cong \pi(\chi'_v)$ if and only if $\chi_v = w \cdot \chi'_v$ for some element $w$ in the Weyl group $W_v = N_{G_v}(T_v)/T_v$. 
\vskip 5pt

(ii)  Every $K_v$-unramified representation of $G(k_v)$  is isomorphic to $\pi(\chi_v)$ for some unramified $\chi_v$.
\vskip 5pt

(iii) Hence the above construction gives a bijection
\[  \{ \text{$K_v$-unramified irreps. of $G(k_v)$} \} \longleftrightarrow \{ \text{unramified characters of $T(k_v)$} \}/ W_v. \]
  \end{prop}   
\vskip 5pt
The proof of this proposition is through showing the Satake isomorphism; the reader can consult \cite[\S 4.22-23]{G1}.
There is an elegant way of reformulating the above proposition, using the Langlands dual group $G_v^{\vee}$  (for split $G_v$) or the Langlands L-group ${^L}G_v$ in general. 
 This reformulation (for split $G_v$ for simplicity) gives a bijection:
 \[  \{ \text{$K_v$-unramified irreps. of $G(k_v)$} \} \longleftrightarrow \{  \text{semisimple conjugacy classes  in $G_v^{\vee}$} \}. \]
 \vskip 10pt
 
\subsection{\bf Tempered representations}
A representation $\pi_v$ is {\em tempered} if its matrix coefficients (which is a function on $G_v$) lies in $L^{2+\epsilon}(G_v)$ for all $\epsilon >0$, i.e. if its matrix coefficients decay sufficiently quickly  \cite[\S 4.27]{G1}. For unramified representations, we can make do with the following ad-hoc definition.
\vskip 5pt
\noindent{\bf Definition:}
A $K_v$-unramified representation $\pi(\chi_v)$ as above is said to be {\em tempered} if  $\chi_v$ is a unitary character, i.e. $|\chi_v| = 1$.
\vskip 5pt
 
 As an example, the trivial representation of $G_v$ is certainly $K_v$-unramified and is contained in the principal series $I(\delta_{B_v}^{1/2})$. Since $\delta_{B_v}^{1/2}$ is not a unitary character, the trivial representation of $G_v$ is not tempered (unless $G_v$ is compact). From the point of view of matrix coefficients, the trivial representation has constant matrix coefficients which certainly do not decay at all.
\vskip 5pt

\subsection{\bf Representation theoretic formulation of the Ramanujan conjecture}

We can now formulate the Ramanujan conjecture for cuspidal representations of a quasi-split group $G$.

\vskip 5pt

\begin{con}[Naive Ramanujan Conjecture] 
Let $\pi = \otimes'_v \pi_v$ be a cuspidal representation of  a quasi-split $G$.  Then for almost all $v$, $\pi_v$ is tempered.
\end{con}

The transition from Ramanujan's original conjecture to this representation theoretic formulation is not clear at all and was first realized by Satake. This transition is described in \cite[\S 4.28]{G1}. The condition that $G$ be quasi-split is there because the conjecture may be easily shown to be false without it. For example, if $G$ is an anisotropic group, so that $[G]$ is compact, then the constant functions, which afford the trivial representation, are certainly cusp forms but the trivial representation is not tempered (as we have remarked above).

\vskip 5pt

\subsection{\bf Counterexample of Howe-PS}
The naive Ramanujan conjecture is expected to be true when $G = \GL_n$, where it has in fact been shown in many cases. 
However, in the 1977 Corvallis proceedings,  Howe and Piatetski-Shapiro \cite{HPS} constructed an example of a cuspidal automorphic representation of the split group $\Sp_4$
which violates the naive Ramanujan conjecture above.   This has led to the following tweak:
\vskip 5pt

\begin{con}[Revised Ramanujan Conjecture]
Let $\pi = \otimes'_v \pi_v$ be a globally generic cuspidal representation of a quasi-split $G$.  Then for almost all $v$, $\pi_v$ is tempered.
\end{con}
Note that for $G = \GL_n$, all cuspidal representations are known to be globally generic (see \cite[\S 6]{G1}).
\vskip 10pt

In this series of lectures, we will follow the same basic idea of Howe-PS and construct a similar counterexample for a quasi-split unitary group $\U_3$.
\vskip 10pt

\subsection{\bf $\epsilon$-Hermitian spaces and unitary groups}
Let us first recall some basics about unitary groups and the underlying Hermitian spaces.
\vskip 5pt

We first begin with an arbitrary field $F$ of characteristic $0$, and let $E$ be an \'etale quadratic $F$-algebra (so $E$ is either a quadratic field extension or  $E= F \times F$), with ${\rm Aut}(E/F)  = \langle c \rangle$ acting on $E$ by $x \mapsto x^c$. With $\epsilon = \pm$, let $V$ be a finite-dimensional  $\epsilon$-Hermitian space over $E$. This means that $V$ is equipped with a nondegenerate $E$-sesquilinear  form $\langle-,- \rangle$, so that
\[  \langle v_1, v_2 \rangle^c = \epsilon \cdot  \langle v_2, v_1 \rangle \quad \text{and} \quad \langle \lambda v_1, v_2 \rangle = \lambda \cdot \langle v_1, v_2 \rangle \]
for $v_1, v_2 \in V$ and $\lambda \in E$. If $\epsilon = +1$, one gets a Hermitian form; if $\epsilon = -1$, one gets a skew-Hermitian form. 
Observe that if $\delta \in E^{\times}$ is a trace 0 element (to $F$), then multiplication-by-$\delta$ takes an $\epsilon$-Hermitian form to a $-\epsilon$-Hermitian form. 
\vskip 5pt

If $(V, \langle-, - \rangle)$ is an $\epsilon$-Hermitian space, let $\U(V)$ be its associated isometry group:
\[ \U(V) = \{ g\in \GL(V): \langle gv_1, gv_2 \rangle = \langle v_1, v_2 \rangle \, \text{for all $v_1, v_2 \in V$} \}. \]
 Because  
\[  \U( V, \langle-, -\rangle) = \U( V, \delta \cdot \langle-, - \rangle) \quad \text{for trace $0$ elements $\delta \in E^{\times}$}, \]
the class of isometry groups one obtains for Hermitian and skew-Hermitian spaces is the same.
These isometry groups are called the {\em unitary groups}: each of them is a connected reductive group with a 1-dimensional anisotropic center $Z \cong E^1$ (the torus defined by the norm 1 elements of $E^{\times}$). 
\vskip 5pt

If $n =\dim_E V$, then  we say that $V$ is maximally split (or simply split) if $V$ contains a maximal isotropic subspace of dimension $[n/2]$. In that case, $\U(V)$ is quasi-split and thus possesses a Borel subgroup $B$ defined over $F$. Such a Borel subgroup is obtained as the stabilizer of a maximal flag of isotropic subspaces:
\[  0 \subset X_1 \subset...\subset X_{[n/2]} \]
with $\dim X_j = j$.

\vskip 5pt

Let us highlight the special case when $E = F \times F$ is the split quadratic $F$-algebra. Then $V$ is an $F\times F$-module and hence has the form $V_0 \times V_0^{\vee}$ for an $F$-vector space $V_0$. 
Up to isomorphism, any Hermitian $E$-space is isomorphic to the one defined by
\[  \langle (v_1, l_1) , (v_2, l_2) \rangle =( l_2(v_1) ,  l_1(v_2)) \in E. \]
Then we note that
 \[  \U(V) \cong \GL(V_0),\] 
via the natural action of $\GL(V_0)$ on $V_0 \times V_0^{\vee}$. We will largely ignore such split cases in the following, as they can be easily handled. 
\vskip 5pt

 \subsection{\bf Invariants of spaces}
 Assume now that $F$ is a local field and  $E/F$ a quadratic field extension.   
A Hermitian space  $V$ of dimension $n$ has  a natural invariant known as the discriminant:
\[
{\rm disc} (V) \in  F^{\times} / \mathrm{N}_{E/F}(E^{\times}). \]
More precisely, if $\{v_1,...,v_n\}$ is an $E$-basis and $A = \left( \langle v_i, v_j \rangle \right)$ is the matrix of inner products of basis elements, then 
\[  {\rm disc}(V) = (-1)^{n(n-1)/2} \cdot \det(A) \in F^{\times}/ NE^{\times}. \]
Using the nontrivial quadratic character $\omega_{E/F}$  of $F^{\times}/N_{E/F}(E^{\times})$, it is convenient to encode $\disc(V)$ as a sign:
\[
\epsilon(V) =  \omega_{E/F}({\rm disc} (V)) = \pm.  \]
 When  $F$ is nonarchimedean, Hermitian spaces are classified by the two invariants $\dim(V)$ and $\epsilon(V)$, so that for each given dimension, there are 2 Hermitian spaces $V^+$ and $V^-$.
When $F = \R$, however, Hermitian spaces are classified by their signatures $(p,q)$, with $p+q = \dim V$.

\vskip 5pt

Likewise, if $W$ is a skew-Hermitian space, then 
\[
 {\rm disc}(W)  \in  \delta^{\dim W} \cdot F^{\times} / \mathrm{N}_{E/F}(E^{\times})
 \]
 and one sets 
 \[
 \epsilon(W)=   \omega_{E/F}(\delta^{-\dim W} \cdot {\rm disc} (W)) =\pm. \]
  Note however that $\epsilon(W)$ depends on the choice of $\delta$.
   \vskip 10pt

Assume now that $F = k$ is a number field and $E/k$ is a quadratic field extension.  Then a Hermitian space $V$ over $E$ is determined by its localizations $\{ V \otimes_k k_v \}$ as $v$ runs over all places of $k$; in other words,  the Hasse principle holds.  Note that  half the places $v$ will split in $E$ and for these, the local situation is the split case (so the split case cannot be ignored for global considerations).
A family of local Hermitian space $\{V_v\}$ (relative to $E_v/k_v$)  arises as the family of localizations of a global Hermitian space  relative to $E/k$ if and only if:
\begin{itemize}
\item for almost all $v$, $\epsilon(V_v) = +$;
\item $\prod_v \epsilon(V_v) = 1$. 
\end{itemize}
There is an analogous statement for skew-Hermitian spaces  which can be formulated, using the observation that multiplication by a nonzero trace 0 element of $E$ switches Hermitian spaces and skew-Hermitian spaces.
 
\vskip 5pt

\subsection{\bf Examples}
Let us consider some examples in small dimension relative to a quadratic field extension of local fields $E/F$. 
\vskip 5pt

\begin{itemize}
\item When $\dim V = 1$, one may identify  $V$ with $E$ (by the choice of a basis element), and  a Hermitian form is given by $(x,y) \mapsto a xy^c$, with $a \in F^{\times}$; 
we denote this rank 1 Hermitian space by $\langle a \rangle$. Then $V^+ = \langle 1 \rangle $ and  $\langle a \rangle \cong \langle b \rangle$ if and only if $a/b \in NE^{\times}$. In any case, $\U( \langle a \rangle ) = E^1 \subset E^{\times}$ for any $a\in F^{\times}$.
\vskip 5pt

We take this occasion to take note of a canonical isomorphism (given by Hilbert's Theorem 90):
\[  E^{\times}/ F^{\times} \cong E^1 \]
defined by
\[   x \mapsto x/x^c. \]
We will frequently use this isomorphism to identify $E^1$ as $E^{\times}/F^{\times}$. 

 \vskip 5pt

\item When $\dim V = 2$, $V^+$ is the split Hermitian space $V^+ = E e_1 \oplus E e_2$, such that
\[  \langle e_1,e_1 \rangle = \langle e_2, e_2 \rangle = 0\quad \text{and} \quad \langle e_1, e_2 \rangle = 1. \]
This 2-dimensional split Hermitian space is also called a hyperbolic plane and is sometimes denoted by $\mathbb{H}$. 
The stabilizer of the isotropic line $E e_1$ is a Borel subgroup, containing a maximal torus  
\[ T = \{  t(a)  = \left( \begin{array}{cc}
a & \\
 & (a^c)^{-1} \end{array} \right): \, a \in E^{\times} \}  \]
and with unipotent radical 
\[  U = \{ u(z)=  \left( \begin{array}{cc}
1 &z \\
 & 1 \end{array} \right) : z \in E,\, Tr_{E/F}(z) = 0 \}. \]
  \vskip 5pt

The other Hermitian space $V^-$ is anisotropic (it has no nonzero isotropic vector). One can describe it in terms of the unique quaternion division $F$-algebra $D$. 
For this, one fixes an $F$-algebra embedding $E \hookrightarrow D$ (unique up to conjugacy by $D^{\times}$ by the Skolem-Noether theorem) and regard $D$ as a 2-dimensional $E$-vector space by left multiplication. One can find $d \in D$ such that $d$ normalizes  $E$ and $dxd^{-1} = x^c$ for $x \in E$, and write $D = E \cdot 1 \oplus E \cdot d$. 
Then one defines a Hermitian form on $D$ by
\[  \langle x,y \rangle = \text{projection of $x \cdot y$ onto $E \cdot 1$.} \]
In terms of this model, one can describe the unitary group $\U(V^-)$ as :
\[  U(V^-) \cong (E^{\times} \times D^{\times})^1/ \nabla F^{\times} = \{ (e, d):  N_E(e) \cdot N_B(b) =1 \}/ \nabla F^{\times}, \]
where $\nabla(F^{\times}) = \{ (t, t^{-1}): t \in F^{\times} \}$. The element $(e,d) \in E^{\times} \times B^{\times}$ gives rise to the operator 
\[  x \mapsto e \cdot x \cdot b^{-1} \]
on $D$. 
\vskip 5pt

 If one replaces $D$ by the matrix algebra $M_2(F)$, the above description of $V^-$ and its isometry group $\U(V^-)$ gives rise to a description of $V^+ = \mathbb{H}$ and $\U(V^+)$. 
This shows that $\U(V^+)$ is intimately connected with $\GL_2$ and $\U(V^-)$ with $D^{\times}$.
\vskip 5pt

\item Consider now the case when $\dim V = 3$. Let $\mathbb{H}$ denote the hyperbolic plane introduced above and recall the 1-dimensional Hermitian space $\langle a \rangle$. 
The sum $\langle a \rangle \oplus \mathbb{H}$ is then a 3-dimensional Hermitian space with $\epsilon( \langle a\rangle + \mathbb{H}) = \omega_{E/F}(a)$.
 As $a$ runs over $F^{\times}/NE^{\times}$, one obtains the two equivalence classes of 3-dimensional Hermitian spaces (in the nonarchimedean case). Thus, both these spaces are split and since $V^+ \cong a \cdot V^-$ for $a \in F^{\times} \setminus NE^{\times}$, we see that $\U(V^+) \cong  \U(V^-)$ is a quasisplit group.  \vskip 5pt

Since this unitary group will play a big role in the Howe-PS construction, let us set up some more notation about it.  Let $\langle a \rangle = E \cdot  v_0$ and $\mathbb{H} = E e  \oplus E e^*$ with $e$ and $e^*$ isotropic vectors. Then with respect to the basis  $\{ e, v_0, e^*\}$ of $\langle a \rangle \oplus \mathbb{H}$, the inner product matrix takes the form
\[  \left( \begin{array}{ccc}
 &  &  1 \\
 & a &  \\
 1 & &    \end{array}  \right). \]
 The Borel subgroup $B = TU$ stabilizing the isotropic line $E \cdot e$ is then upper triangular, with elements of $T$ and $U$ taking the form
 \[  t(a, b)  = \left( \begin{array}{ccc}
 a & & \\
  & b &  \\
   &  & (a^c)^{-1}  \end{array} \right) \quad \text{with $a \in E^{\times}$ and $b \in E^1$}  \]
and
\[  u(x,z) = \left( \begin{array}{ccc}
1 & 0& z \\
& 1 & 0 \\
&  &  1  \end{array}  \right) \cdot 
 \left( \begin{array}{ccc}
1 & -a^c x^c   &  *  \\
& 1 & x \\
&  &  1  \end{array}  \right), \]
 with $x ,z \in E$ and $Tr_{E/F}(z) = 0$.
\end{itemize}

\vskip 5pt

\subsection{\bf Basic idea of Howe-PS construction}
We can now give a brief summary of the Howe-PS construction of cuspidal representations of $\U_3$ which violate the naive Ramanujan conjecture. 
\vskip 5pt

Let $E/k$ be a quadratic field extension of number fields. 
Consider a skew-Hermitian space $(W, \langle -, - \rangle)$ of dimension $3$ over $E$.  We would like to produce some cusp forms on $\U(W)$ which violates the naive Ramanujan conjecture.  These functions on $\U(W)$ will be obtained by restriction (or pullback) of a simpler class of automorphic forms on a larger group containing $\U(W)$. What is this larger group?  
\vskip 5pt

By restriction of scalars, we have a 6-dimensional space ${\rm Res}_{E/k}(W)$ over $k$ and the $k$-valued form ${\rm Tr}_{E/k}(\langle-, -\rangle)$ defines a symplectic form on ${\rm Res}_{E/k}(W)$ with associated symplectic group $\Sp({\rm Res}_{E/k}(W))$. This defines an embedding
\[  \iota: \U(W) \hookrightarrow \Sp({\rm Res}_{E/k}(W)). \]
It turns out that the simple automorphic forms we need are not really living on $\Sp({\rm Res}_{E/k}(W))$.
Rather,  the symplectic group $\Sp_{2n}(\A)$ has a topological $S^1$-cover $\Mp_{2n}(\A)$ known as the metaplectic group (where $S^1$ is the unit circle in $\C^{\times}$):
\[  \begin{CD}
1 @>>> S^1 @>>> \Mp_{2n}(\A) @>>>\Sp_{2n}(\A) @>>>1 \end{CD} \]
and the simpler automorphic forms in question actually live on $\Mp_{2n}(\A)$. The need to work with this nonlinear cover accounts for much of the technicality of this subject, but one cannot argue with nature.

 \vskip 5pt

These simpler automorphic forms on the metaplectic groups are the theta functions and the automorphic representations they span are called the Weil representations.
In order to pull back these theta functions from $\Mp({\rm Res}_{E/k}(W))(\A)$ to $\Sp({\rm Res}_{E/k}(W))$,  one needs to construct a lifting of $\iota$ to
\[ \tilde{\iota}: \U(W)(\A) \hookrightarrow \Mp({\rm Res}_{E/k}(W))(\A). \]
This is highly technical but it can be done and Howe-PS then restricted these theta functions to $\U(W)(\A)$.
\vskip 5pt

Representation theoretically, if 
\[  \Omega \subset  \mathcal{A}(\Mp({\rm Res}_{E/k}(W))) \]
denotes an automorphic Weil representation, and 
\[  \tilde{\iota}^* : \mathcal{A}(\Mp({\rm Res}_{E/k}(W))) \rightarrow \mathcal{A}(\U(W)) \]
denotes the restriction of functions, then one obtains a $\U(W)(\A)$-submodule
\[  \tilde{\iota}^*(\Omega) \subset \mathcal{A}(\U(W)). \]
Now recall that the center of $\U(W)$ is isomorphic to $E^1$ as an algebraic group. One can spectrally decompose $\tilde{\iota}^*(\Omega)$  according to central characters.
\[  \tilde{\iota}^*(\Omega) =  \bigoplus_{\chi} \Omega_{\chi} \]
as $\chi$ runs over the automorphic characters of $E^1$, or equivalently of $E^{\times}/k^{\times}$, where 
\[  \Omega_{\chi} = \{ f \in \Omega: f(zg) = \chi(z) \cdot f(g) \, \text{for all $z \in Z(\U(W)) =\A_E^1$ and $g \in \U(W)(\A)$} \}. \]
What we would like to show is that:
\vskip 5pt

\begin{itemize}
\item  for each $\chi$, $\Omega_{\chi}$ is  an irreducible automorphic representation of $\U(W)$ and is cuspidal for many $\chi$.
\item for any $\chi$, $\Omega_{\chi}$   violates the naive Ramanujan conjecture.
\end{itemize}
One can view the map $\chi \mapsto \Omega_{\chi}$ as a lifting of automorphic representations from $\U_1$ to $\U_3$. This lifting is an instance of the theta correspondence, which we will discuss in the next two lectures.

\vskip 10pt

\newpage

\section{\bf Lecture 2: Local Theta Correspondence}
The next two lectures will be devoted to a discussion of the theory of theta correspondence, so as to understand the construction of Howe-PS in its proper context.
Two possible references for this are the survey papers of D. Prasad \cite{P} and S. Gelbart \cite{Ge}. In particular, the latter is concerned with theta correspondence for unitary groups.
 In this second lecture, we shall focus on the local theta correspondence, for which a basic reference is the book \cite{MVW} of Moeglin-Vigneras-Waldspurger.
 Hence, we will be working over a local field $F$ (of characteristic 0), and for simplicity, we shall assume $F$ is nonarchimedean.
\vskip 5pt

\subsection{\bf Basic idea}
From the last lecture, we saw that the Howe-PS  construction basically gives a map 
\[  \{ \text{Automorphic characters of $\U_1 = E^1$} \} \rightarrow \{ \text{Automorphic representations of $\U_3$} \} \]
 sending $\chi$ to $\Omega_{\chi}$.  
 \vskip 5pt
 
 Now given any two groups $G$ and $H$, one may ask more generally: what are some ways of constructing such a lifting from ${\rm Irr}(G)$ to ${\rm Irr}(H)$? Here, ${\rm Irr}(G)$ denotes the set of equivalence classes of irreducible representations of $G$.
\vskip 5pt

A standard procedure is as follows. Suppose for simplicity that $G$ and $H$ are finite groups and  $\Omega$ is a (finite-dim) representation of $G \times H$. Then one may decompose $\Omega$ into irreducible $G \times H$-summands:
\[  \Omega =  \bigoplus_{\pi \in {\rm Irr}(G)}  \bigoplus_{\sigma \in {\rm Irr}(H)}   m(\pi, \sigma) \cdot \pi \otimes \sigma. \]
One can rewrite this as:
\[  \Omega = \bigoplus_{\pi \in {\rm Irr}(G)} \pi \otimes V(\pi) \]
where 
\[  V(\pi) =  \bigoplus_{\sigma \in {\rm Irr}(H)} m(\pi, \sigma) \sigma. \] 
This gives a map  
\[  {\rm Irr}(G)  \longrightarrow R(H) \, \, \text{ (Grothendieck group of ${\rm Irr}(H)$}) \]
sending $\pi$ to $V(\pi)$.  Since we are interested in getting irreducible representations of $H$ as outputs, we ask: for what $\Omega$ is $V(\pi)$ is  irreducible or zero for any $\pi$? 
If we can find such an $\Omega$, then the map $\pi \mapsto V(\pi)$ would be a map
\[  {\rm Irr}(G) \longrightarrow {\rm Irr}(H) \cup \{ 0\}. \]
\vskip 5pt

An example of an $\Omega$ that has this property is certainly the trivial representation. However the map so obtained is not very interesting. On the other hand, if $\dim \Omega$ is too big, then $\dim V(\pi)$ will have to be big for many $\pi$'s as well, so that $V(\pi)$ tends to be reducible. Thus, a simple heuristic is that $\Omega$ cannot be too big nor too small. 

\vskip 5pt

In practice, one can try an $\Omega$ arising in the following way.  Suppose there is an (almost injective)  group homomorphism 
 \[ \iota: G \times H \longrightarrow E \quad \text{for some group $E$.} \]
 One can take $\Omega$ to be an irreducible representation of $E$ of smallest possible dimension $>1$ and then pull it back to $G \times H$. 
\vskip 5pt

The theory of theta correspondence, which was systematically developed by Howe, arises in this way.   
\vskip 10pt

\subsection{\bf Reductive dual pairs}
 Let  $F$ be a field of characteristic $0$, and let $E$ be an \'etale quadratic $F$-algebra, with ${\rm Aut}(E/F)  = \langle c \rangle$.  Let $V$ be a finite-dimensional Hermitian space
over $E$ and $W$ a skew-Hermitian space. Then $V \otimes_E W$ is naturally a skew-Hermitian space over $E$. By restriction of scalars, we may regard $V \otimes_E W$ as an $F$-vector space which is equipped with  a natural symplectic form ${\rm Tr}(\langle -, - \rangle_V \otimes \langle - , - \rangle_W)$.  
Then one has a natural map of isometry groups
\[ \iota:  {\rm U}(V)  \times {\rm U}(W)  \longrightarrow {\rm Sp}(V \otimes_E W) \]
which is injective on each factor $\U(V)$ and $\U(W)$.
The images of ${\rm U}(V)$ and ${\rm U}(W)$ are in fact mutual commutants of each other in the symplectic group, and such a pair of groups is called a {\em reductive dual pair}. 
\vskip 5pt
Howe has given a complete classification of all such dual pairs in the symplectic group. A further example (perhaps easier than the one above) is obtained as follows.
If $V$ is a symmetric bilinear space (or a quadratic space) and $W$ a symplectic space over $F$, then $V \otimes_F W$ inherits a natural symplectic form (by tensor product) and one has
\[  \O(V) \times \Sp(W) \longrightarrow \Sp(V \otimes W). \]

\vskip 10pt

\subsection{\bf Metaplectic groups and Heisenberg-Weil representations}

Assume now that $F$ is a local field.  The symplectic group $ {\rm Sp}(V \otimes_E W)$ has a nonlinear $S^1$-cover ${\rm Mp}(V \otimes_E W)$ known as the metaplectic group:
\[  \begin{CD}
 1 @>>> S^1 @>>> \Mp(V \otimes_E W) @>>> \Sp(V \otimes_E W) @>>> 1 \end{CD} \]
The construction of this central extension is a lecture course in itself, but since its construction is such a classic result, we feel obliged to give a sketch. 
\vskip 5pt

Let us work in the context of an arbitrary symplectic vector space $W$ over $F$  (in place of the cumbersome notation $V \otimes_E W$). Let
\[  H(W)  = W \oplus F \]
and equip it with the group law
\[  (w_1, t_1) \cdot (w_2, t_2) = (w_1 + w_2, t_1 + t_2 + \frac{1}{2} \langle w_1, w_2\rangle). \]
This is the so-called Heisenberg group. It is a 2-step nilpotent group with center $Z = F$ and  commutator $[H(W), H(W)] = Z$. 
The definition of this group law is motivated by the Heisenberg commutator relations from quantum mechanics, hence the name. 
\vskip 5pt

The irreducible (smooth) representations of $H(W)$ can be classified and come in 2 types:
\vskip 5pt

\begin{itemize}
\item the 1-dimensional representations: these factor through 
\[  H(W) / [H(W), H(W)] = H(W)/Z   \cong W \]
 and hence are given by characters of $W$. Observe that these are precisely the representations trivial on the center $Z$.
\vskip 5pt

\item  the other irreducible representations have nontrivial central character. For a fixed $\psi: Z = F \longrightarrow \C^{\times}$, the Stone-von Neumann theorem asserts that $H(W)$ has a unique irreducible representation $\omega_{\psi}$ with central character $\psi$. Moreover, the representation $\omega_{\psi}$ is unitary. 
\end{itemize}

One can give an explicit construction of $\omega_{\psi}$. Let 
\[  W = X  \oplus Y \]
be a Witt decomposition of $W$ so that $X$ and $Y$ are maximal isotropic subspaces.  Then $H(X) = X + F$ is an abelian subgroup of $H(W)$ and one can extend $\psi$ to a character of $H(X)$ by setting $\psi(x,0) = 1$ (i.e. extending trivially to $X$).  Then
\[   \omega_{\psi} \cong {\rm ind}_{H(X)}^{H(W)} \psi  \quad \text{(compact induction).} \]
This induced representation is realized on $\mathcal{S}(Y) := C^{\infty}_c(Y)$, and the action of $h \in H(W)$ can be easily written down:
\[ \begin{cases}
(\omega_{\psi}(0,z) f) (y') = \psi(z) \cdot f(y'), \text{  for $z \in F$;} \\
(\omega_{\psi}(y,0) f)(y') = f(y+y'),  \text{  for $y \in Y$} \\
(\omega_{\psi}(x,0) f)(y') = \psi(\langle y',x\rangle) \cdot f(y'),  \text{  for $x \in X$.}
\end{cases} \]
This action preserves the natural inner product on $\mathcal{S}(Y)$.
\vskip 5pt

The symplectic group $\Sp(W)$ acts on $H(W)$ as group automorphisms via:
\[  g \cdot (w, t) = (g \cdot w, t). \]
Observe that the action on $Z$ is trivial.  Hence the representation ${^g}\omega_{\psi} = \omega_{\psi} \circ g^{-1}$ is irreducible and has the same nontrivial central character as $\omega_{\psi}$. By the Stone-von Neumann theorem, these two representations are isomorphic, i..e there exists an invertible operator  $A_{\psi}(g)$ on the underlying vector space $\mathcal{S}$  of $\omega_{\psi}$ such that
\[    A_{\psi}(g) \circ {^g}\omega_{\psi}(h)    =  \omega_{\psi}  \circ A_{\psi}(g) \quad \text{for all $h \in  H(W)$.} \]
By Schur's lemma, the operator $A_{\psi}(g)$ is well-defined up to $\C^{\times}$. By the unitarity of $\omega_{\psi}$, we can insist that $A_{\psi}(g)$ is unitary and hence it is well-defined up to the unit circle $S^1 \subset \C^{\times}$.  Hence we have a map
\[  A_{\psi}: \Sp(W) \longrightarrow \GL(\mathcal{S}) / S^1. \] 
When one pulls back the short exact sequence
\[  \begin{CD}
1 @>>> S^1 @>>> \GL(\mathcal{S}) @>>>  \GL(\mathcal{S}) /S^1 @>>> 1 \end{CD} \]
by the map $A_{\psi}$, one obtains the desired metaplectic group $\Mp_{\psi}(W)$:
\[  \begin{CD}
1 @>>> S^1 @>>>  \Mp_{\psi}(W) @>>> \Sp(W) @>>> 1 \\
 @. @| @VV{\tilde{A}_{\psi}}V @VV{A_{\psi}}V  \\
 1 @>>> S^1 @>>> \GL(\mathcal{S}) @>>> \GL(\mathcal{S})/ S^1 @>>> 1 
 \end{CD} \]
 Hence, this construction produces not just the group $\Mp_{\psi}(W)$ but also a natural representation
\[  \tilde{A}_{\psi}:  \Mp_{\psi}(W) \longrightarrow  \GL(\mathcal{S}). \]
Thus, we have a representation $\omega_{\psi}$ of $\Mp_{\psi}(W) \ltimes H(W)$. In other words, the irreducible representation $\omega_{\psi}$ of $H(W)$ extends to the semidirect product $H(W) \rtimes \Mp_{\psi}(W)$ (with $\Mp_{\psi}(W)$ acting on $H(W)$ via its projection to $\Sp(W)$).  We call this a Heisenberg-Weil representation; its restriction to $\Mp_{\psi}(W)$ is simply called a Weil representation. 
\vskip 10pt

It turns out that the isomorphism class of the extension defining $\Mp_{\psi}(W)$ is independent of $\psi$; so we shall write $\Mp(W)$ henceforth, suppressing $\psi$.
The Weil representation $\omega_{\psi}$ of $\Mp(W)$  is, in some sense, the smallest infinite-dimensional representation of the metaplectic group. We only make two further remarks here:

\vskip 5pt

\begin{itemize}
\item $\omega_{\psi}$ is a genuine representation of $\Mp(W)$, in the sense that $\omega_{\psi}(z) = z$ for all $z \in S^1$, so that $\omega_{\psi}$ does not factor to a smaller cover of $\Sp(W)$. 

\item $\omega_{\psi}$ is not an irreducible representation of $\Mp(W)$, but rather a direct sum of two irreducible representations (the even and odd Weil representaitons):
\[  \omega_{\psi} = \omega_{\psi}^+ \oplus \omega_{\psi}^-. \]
Indeed, with $\omega_{\psi}$ realized on $\mathcal{S}(Y)$ as above, $\omega_{\psi}^+$ is realized on the subspace of even functions  (i.e. $f(-x) = f(x)$) whereas $\omega_{\psi}^-$ is realized on the subspace of odd functions.
\end{itemize}
\vskip 5pt

\subsection{\bf Schrodinger model}
One can ask if it is possible to write down some formulas for the action of elements of $\Mp(W)$, for example on the model $\mathcal{S}(Y)$ of $\omega_{\psi}$. 
Let $P(X)$ be the maximal parabolic subgroup of $\Sp(W)$ stabilizing the maximal isotropic subspace $X$; this is the so-called Siegel parabolic subgroup.
Its Levi decomposition has the form $P(X) = M(X) \cdot N(X)$, with $M(X) \cong \GL(X)$ and 
\[  N(X)  = \{ n(B): B \in Sym^2 X \subset \Hom(Y,X) \}, \]
 where
\[  n(B) = \left( \begin{array}{cc}
1 & B \\
 0 & 1 \end{array} \right)  \quad \text{(relative to $W = X \oplus Y$).} \]
 The metaplectic $S^1$-cover splits canonically over $N(X)$ and noncanonically over $M(X)$.
Then one can write down formulas for the action of elements lying over $g \in \GL(X)$ and $n(B) \in N(X)$:
\[ \begin{cases}
(\omega_{\psi}(g) f )(y) =   |\det_X(g)|^{1/2} \cdot f(g^{-1} \cdot y)   \\
(\omega_{\psi}(n(B)) f)(y)   = \psi( \frac{1}{2} \cdot \langle By, y \rangle) \cdot f(y). \end{cases} \]
 To describe the action of $\Mp(W)$ (at least projectively), one needs to give the action of an extra Weyl group element 
\[  w = \left( \begin{array}{cc}
0 & 1 \\
 -1 & 0 \end{array} \right)  \]
 which together with $P(X)$ generates $\Sp(W)$. The action of this $w$ is given by Fourier transform  (up to scalars).
 \vskip 5pt
 
 The above gives  the Schrodinger model of the Weil representation (which is related to the Schrodinger description of quantum mechanics).
This concludes our brief sketch of the construction of the metaplectic group and its Weil representations.
\vskip 10pt

\subsection{\bf  Weil representations of unitary groups}
Let us return to the setting of our unitary dual pair
\[  \iota: \U(V) \times \U(W) \longrightarrow \Sp(V \otimes_E W). \]
By the above, the metaplectic group $\Mp(V \otimes_E W)$  has a distinguished representation $\omega_{\psi}$ depending on a nontrivial additive character $\psi$ of $F$.  
 If the embedding $\iota$ can be lifted to a homomorphism 
\[  \tilde{\iota}:  {\rm U}(V)  \times  {\rm U}(W)  \longrightarrow {\rm Mp}(V \otimes_E W), \]
then we obtain a representation $\omega_{\psi} \circ  \tilde{\iota}$ of ${\rm U}(V) \times {\rm U}(W)$.  
 
 \vskip 5pt
 
  Such splittings have been constructed and classified by S. Kudla \cite{K}.    They are not unique but 
 can be specified  by picking two characters $\chi_V$ and $\chi_W$ of $E^{\times}$ such that
\[  \chi_V |_{F^{\times}} =  \omega_{E/F}^{\dim V} \quad \text{and}  \quad \chi_W|_{F^{\times}}  = \omega_{E/F}^{\dim W}. \]
One way of doing this is, for example, fixing a character $\gamma$ of $E^{\times}$ such that $\gamma|_{F^{\times}} = \omega_{E/F}$ (i.e. a conjugate symplectic character) and then taking
\[  \chi_V = \gamma^{\dim V} \quad \text{and} \quad \chi_W = \gamma^{\dim W}. \]
In any case, 
for a fixed  pair $(\chi_V, \chi_W)$ of splitting characters,  Kudla provides a splitting
\[  \tilde{\iota}_{\chi_V, \chi_W, \psi} : \U(V) \times \U(W) \longrightarrow  \Mp(V \otimes_E W) \]
of $\iota$. In fact, the choice of $\chi_V$ gives rise to a splitting
\[  \iota_{V,W, \chi_V, \psi}: \U(W) \longrightarrow \Mp(V \otimes_E W) \]
over $\U(W)$, whereas the choice of $\chi_W$ gives a splitting
\[  \iota_{V,W, \chi_W, \psi} :\U(V) \longrightarrow \Mp(V \otimes_E W). \]
Hence, the splitting over the two members of the dual pair can be constructed somewhat independently of each other.  This is a manifestation of a basic property of the metaplectic cover: if two elements of $\Sp(V \otimes_E W)$ commute, then any lifts of them in $\Mp(V \otimes_E W)$ also commute with each other. 
\vskip 5pt

With a splitting fixed,  we may pull back the Weil representation $\omega_{\psi}$ and obtain a representation
\[  \Omega_{V, W, \chi_V, \chi_W, \psi} := \tilde{\iota}_{\chi_V, \chi_W, \psi}^*(\omega_{\psi}) \quad \text{ of $\U(V) \times \U(W)$.}\]
We call this $\Omega_{\chi_V, \chi_W, \psi}$ a Weil representation of $\U(V) \times \U(W)$; we will often suppress the subscript for ease of notation.
 
\vskip 5pt
While we have not described the splitting $\tilde{\iota}_{\chi_V, \chi_W, \psi}$ explicitly, we highlight a few basic properties:
\vskip 5pt

\begin{itemize}
\item (Change of $(\chi_V, \chi_W)$) One can easily describe the effect of changing $(\chi_V, \chi_W)$ on $\Omega_{\chi_V, \chi_W, \psi}$. More precisely, if $(\chi_v' , \chi_W')$ is another pair of splitting characters, then
\[  \Omega_{V,W,\chi_V', \chi_W', \psi} \cong   \Omega_{V,W, \chi_V, \,\chi_W, \psi}  \otimes (\chi_V'/\chi_V \circ i \circ {\det}_W)  \otimes (\chi_W'/\chi_W \circ i \circ {\det}_V), \]
where $i: E^1 \rightarrow E^{\times}/F^{\times}$, taking note that $\chi_V'/\chi_V$ and $\chi_W'/\chi_W$  are characters of $E^{\times}/F^{\times}$.
 \vskip 5pt
 
 \item (Scaling)  Given $a \in F^{\times}$, one can scale the additive character $\psi$ to obtain  $\psi_a(x) = \psi(ax)$. One can also scale the Hermitian or skew-Hermitian forms on $V$ and $W$, otaining $V^a$ and $W^a$. If we identify $\U(V)$ and $\U(V^a)$ as the same subgroup of $\GL(V)$, then one has: 
 \[  \Omega_{V, W, \chi_V, \chi_W, \psi_a} \cong \Omega_{V^a, W, \chi_V, \chi_W, \psi} \cong \Omega_{V, W^a, \chi_V, \chi_W, \psi}. \]
 \vskip 5pt
 
 \item (Duality) Recalling that the Weil representation is unitarizable, so that its dual is isomorphic to its complex conjugate, we note:
 \[ \overline{ \Omega_{V,W,\chi_V, \chi_W, \psi}} \cong \Omega_{V,W, \chi_V^{-1}, \chi_W^{-1}, \overline{\psi}}. \]
 \vskip 5pt
 
 \item (Center) If we take 
  \[  \chi_V = \gamma^{\dim V} \quad \text{and} \quad \chi_W = \gamma^{\dim W} \]
  for a conjugate symplectic character $\gamma$, and identify $Z_{\U(V)}$ and $Z_{\U(W)}$ with $E^1$ (and hence with each other), 
  $ \tilde{\iota}_{\chi_V, \chi_W, \psi}$ agrees on $Z_{\U(V)}$ and $Z_{\U(W)}$. 
  \end{itemize}
  \vskip 5pt

\vskip 5pt
\subsection{\bf Local theta lifts}
In the theory of local theta correspondence, one studies  how the representation $\Omega_{V, W, \chi_V, \chi_W, \psi}$ decomposes into irreducible pieces. 
For this, one would need to know more about the representation $\Omega_{V,W, \chi_V, \chi_W, \psi}$ than what has been described above. For example, one may demand if there are formulas for the group action. Ultimately, such formulas would have been derived from those in the Schrodinger model we described earlier and an explicit knowledge of the splitting $\tilde{\iota}_{\chi_V, \chi_W, \psi}$.  We will come to this in the particular case of interest later on in this lecture. At the moment, let us just formulate some questions one may ask and  describe the answers to some of these questions.

\vskip 10pt

We will write ${\rm Irr}({\rm U}(W))$ for the set of equivalence classes of irreducible smooth representations of ${\rm U}(W)$.
Unlike the case of finite (or compact) groups, the representation $\Omega$ is infinite-dimensional and is not necessarily semisimple as a $\U(V) \times \U(W)$-module. So when one talks about the decomposition of $\Omega$ into irreducible constituents, one can understand this in potentially two ways:
\begin{itemize}
\item for which $\pi \otimes \sigma \in {\rm Irr}(\U(V)) \times {\rm Irr}(\U(W))$ is $\pi \otimes \sigma$ a subrepresentation of $\Omega$? 
\item for which $\pi \otimes \sigma \in {\rm Irr}(\U(V)) \times {\rm Irr}(\U(W))$ is $\pi \otimes \sigma$ a quotient of $\Omega$? 
\end{itemize}
It turns out that it is more fruitful to consider the second question above. Hence, $\Omega$ determines a subset of ${\rm Irr}(\U(V) \times {\rm Irr}(\U(W))$:
\[  \Sigma_{\Omega} = \{ (\pi, \sigma): \text{$\pi \otimes \sigma$ is a quotient of $\Omega$}\}  \]
and thus a correspondence between ${\rm Irr}(\U(V))$ and  ${\rm Irr}(\U(W))$. This is the correspondence in ``theta correspondence". 
If $(\pi , \sigma) \in \Sigma_{\Omega}$, we say that $\sigma$ is a local theta lift of $\pi$ and vice versa. 
\vskip 5pt

One can reformulate the above definition in a slightly different way, which is more convenient for the question of determining all possible theta lifts of a given $\pi$. 
For $\pi \in {\rm Irr}({\rm U}(V))$, one considers the maximal $\pi$-isotypic quotient of $\Omega$:
\[  \Omega /  \bigcap_{f \in \Hom_{\U(V)}(\Omega, \pi)}  {\rm Ker} (f),  \]
which is a $\U(V) \times \U(W)$-quotient of $\Omega$ expressible in the form
\[ \pi \otimes \Theta(\pi), \]
for some smooth representation $\Theta(\pi)$ of $\U(W)$ (possibly zero, and possibly infinite length a priori).
 We call $\Theta(\pi)$ the {\em big theta lift} of $\pi$. An alternative way to define $\Theta(\pi)$ is:
\[  \Theta(\pi) = ( \Omega \otimes \pi^{\vee})_{\U(V)}, \]
 the maximal $\U(V)$-invariant quotient of $\Omega \otimes \pi^{\vee}$.  To see the equivalence of these two definitions of $\Theta(\pi)$:
 \vskip 5pt
 
 \begin{itemize}
 \item we have a natural equivariant projection
 \[  q: \Omega \longrightarrow  \Omega /  \bigcap_{f \in \Hom_{\U(V)}(\Omega, \pi)}  {\rm Ker} (f) = \pi \otimes \Theta(\pi). \]
 In addition, the natural map
 \[  \Omega \otimes \pi^{\vee} \longrightarrow (\Omega \otimes \pi^{\vee})_{\U(V)} \]
 induces another equivariant map
 \[  p: \Omega \longrightarrow \pi \boxtimes   (\Omega \otimes \pi^{\vee})_{\U(V)}.\]
 \vskip 5pt
 \item note that
 \[  \Hom_{\U(V)}(\Omega, \pi) \cong \Hom( (\Omega \otimes \pi^{\vee})_{\U(V)}, \C). \] 
 For $f \in \Hom_{\U(V)}(\Omega, \pi)$, with corresponding $\ell_f \in  (\Omega \otimes \pi^{\vee})^*$, one has a commutative diagram
 \[  \begin{CD}
 \Omega @>p>>  \pi \boxtimes (\Omega \otimes \pi^{\vee})_{\U(V)} \\
 @VfVV   @VV\ell_fV \\
 \pi  @= \pi 
 \end{CD}
\]  
 This implies that 
\[  {\rm Ker}(p) \subset \bigcap_{f \in \Hom_{\U(V)}(\Omega, \pi)}  {\rm Ker} (f)  =  {\rm Ker}(q). \]
 \vskip 5pt
 
 \item conversely, if $\phi \in {\rm Ker}(q)$, then $p(\phi) \in \pi \boxtimes  (\Omega \otimes \pi^{\vee})_{\U(V)}$ is killed by every element of $(\Omega \otimes \pi^{\vee})_{\U(V)}^*$, so that $p(\phi) = 0$. 
 \end{itemize}

 In any case, it follows from definition that there is a natural $\U(V) \times \U(W)$-equivariant map
 \[  \Omega \twoheadrightarrow \pi \otimes \Theta(\pi), \]
 which satisfies the ``universal property" that for any smooth representation $\sigma$ of $\U(W)$, 
 \[  \Hom_{\U(V) \times \U(W)}(\Omega, \pi \otimes \sigma) \cong \Hom_{\U(W)}(\Theta(\pi), \sigma)  \quad \text{(functorially)}. \] 
 The local  theta lifts of $\pi$ are then the irreducible quotients of $\Theta(\pi)$.
 \vskip 5pt

\subsection{\bf Howe duality conjecture}
 
The goal of local theta correspondence is to determine the representation $\Theta(\pi)$ or rather its irreducible quotients.
Recall that our hope is that $\Theta(\pi)$ is close to being  irreducible or at least not too big. To this end, we first note the following basic result of Howe and Kudla:
\begin{prop}[Finiteness]
For any $\pi \in {\rm Irr}(\U(V))$, $\Theta(\pi)$ is of finite length. In particular, if $\Theta(\pi)$ is nonzero, it has (finitely many) irreducible quotients and we may consider its maximal semisimple quotient (its cosocle) $\theta(\pi)$. Moreover, for any $\pi \in{\rm Irr}(\U(V))$ and $\sigma \in {\rm Irr}(\U(W))$,
\[  \dim \Hom_{\U(V) \times \U(W)}(\Omega, \pi \otimes \sigma) < \infty. \]
\end{prop}
We call $\theta(\pi)$ the {\em small theta lift} of $\pi$. The local theta lifts of $\pi$ are precisely the irreducible summands of $\theta(\pi)$.  
 \vskip 5pt

We can now formulate a fundamental theorem \cite{W2, GT}:
\vskip 5pt

\begin{thm}[Howe duality theorem]
(i) If $\Theta(\pi)$ is nonzero, then it has a unique irreducible quotient. In other words,  $\theta(\pi)$ is irreducible or $0$.
\vskip 5pt

(ii) If $\theta(\pi) \cong \theta(\pi') \ne 0$, then $\pi \cong \pi'$.
\vskip 5pt

Hence, one has an map
\[  \theta_{\chi_V, \chi_W, \psi} : {\rm Irr}(\U(V)) \longrightarrow  {\rm Irr}(\U(W)) \cup \{ 0 \} \]
which is injective when restricted to the subset of ${\rm Irr}(\U(V))$ which is not sent to $0$. 
Moreover,
\[  \dim \Hom_{\U(V) \times \U(W)}(\Omega, \pi \otimes \sigma) \leq 1 \]
for any $\pi \in {\rm Irr}(\U(V))$ and $\sigma \in {\rm Irr}(\U(W))$. 
\end{thm}
\vskip 5pt
Another way to formulate the theorem is to note that the subset/correspondence  $\Sigma_{\Omega}$ is the graph of a bijective function between $pr_V(\Sigma_{\Omega}) \subset {\rm Irr}(\U(V))$ and $pr_W(\Sigma_{\Omega}) \subset {\rm Irr}(\U(W))$, where $pr_V$ refers to the projection to ${\rm Irr}(\U(V))$. 
Thus, we see that local theta correspondence is an instance of the Basic Idea highlighted at the beginning of the lecture. 
\vskip 10pt

\subsection{\bf Questions}
After the Howe duality theorem above, we can ask the following questions:
\vskip 5pt

\begin{itemize}
\item[(a)] (Nonvanishing)  For a given $\pi \in {\rm Irr}(\U(V))$, decide if $\theta_{\chi_V, \chi_W, \psi}(\pi)$ is nonzero.
\vskip 5pt

\item[(b)] (Identity) Describe the map $\theta_{\chi_V, \chi_W, \psi}$ explicitly. In other words,  if $\theta_{\chi_V, \chi_W, \psi}(\pi)$ is nonzero, can one describe it in terms of $\pi$ in another way? 
\end{itemize}
Nowadays, both these questions have rather complete answers, but it is too much to describe these answers for this course. Instead, we will highlight some relevant answers for our application.
\vskip 5pt

\subsection{\bf Rallis' tower}
For the nonvanishing question,  Rallis observed that it is fruitful to consider theta correspondence in a family. 
Let $W_0$ be an anisotropic skew-Hermitian space over $E$ (for example a 1-dimensional one), and for $r \ge 0$, let
\[  W_r = W_0 \oplus \mathbb{H}^r \]
where $\mathbb{H}$ is the hyperbolic plane.
The collection $ \{ W_r \, | \, r \ge 0 \}$ is called a Witt tower of spaces. Observe that:
\begin{itemize}
\item $\dim W_r \mod 2$ is independent of $r$;
\item  $\disc(W_r)$ or equivalently the sign $\epsilon(W_r)$ is independent of $r$.
\end{itemize}
Hence, in the nonarchimedean case, there are precisely two Witt towers of skew-Hermitian spaces with a fixed dimension modulo 2, and any
given skew-Hermitian space $W$ is a member of a unique Witt tower. 
\vskip 5pt

One can then consider a family of theta correspondences associated to the tower of reductive dual pairs 
$({\rm U}(V), {\rm U}(W_r))$ with respect to a fixed pair of splitting characters $(\chi_V, \chi_W)$ (note that we can fix $\chi_{W_r}$ independently of $r$, since the parity of $\dim W_r$ is independent of $r$).  Kudla showed \cite{K}:
\vskip 5pt

\begin{prop}
(i) For $\pi \in  {\rm Irr}({\rm U}(V))$, there is a smallest $r_0 = r_0(\pi)$ such that $\Theta_{V, W_{r_0}, \psi}(\pi)\ne 0$. Moreover, $r_0 \leq \dim V$.

(ii) For any $r > r_0$,   $\Theta_{V, W_r, \psi}(\pi)\ne 0$ (tower property).

(iii) Suppose that $\pi$ is supercuspidal. Then $\Theta_{V, W_{r_0}, \psi}(\pi)$ is irreducible supercuspidal. For $r > r_0$, $\Theta_{V, W_r, \psi}(\pi)$ is irreducible but not supercuspidal.
\end{prop}
 Some remarks:
 
 \vskip 5pt
 
 \begin{itemize}
 \item We call this smallest $r_0(\pi)$ from the proposition above the {\em first occurrence index} of $\pi$ in the Witt tower $(W_r)$ (relative to our fixed  $(\chi_V, \chi_W, \psi)$). 
 
 \item The fact that $r_0 \leq  \dim V$ means that when $W$ is sufficiently large, more precisely when $W$ has an isotropic  subspace of dimension $\geq \dim V$, the map $\theta_{V,W, \psi}$ is nonzero on the whole of ${\rm Irr}(\U(V))$. When $r \geq \dim V$, we say that the theta lifting is in the {\em stable range}.

\item The nonvanishing problem (a) highlighted above is reduced to the question of determining the first occurrence indices for the two Witt towers.
\end{itemize}
 \vskip 5pt
In fact, our job is halved because the two first occurrence indices (for a given $\dim W_{\bullet} \mod 2$) are not independent of each other. Rather, we have the following theorem of Sun-Zhu \cite{SZ}:
\vskip 5pt

\begin{thm}[Conservation relation]
Consider the two towers $(W_v)$ and $(W'_r)$ of skew-Hermitian spaces with fixed $\dim W_{\bullet} \mod 2$, and let $r_0$ and $r_0'$ be the respective first occurrence indices of a fixed $\pi \in {\rm Irr}(\U(V))$ (relative to a fixed data $(\chi_V, \chi_W, \psi)$). Then
\[  \dim W_{r_0} + \dim W'_{r'_0} =  2 \dim V + 2. \]
In particular $r_0$ and $r'_0$ determine each other.
\end{thm}

The  conservation relation above implies the following dichotomy theorem (which you should try to prove):

\begin{cor}  \label{C:dichotomy}
Suppose that $W$ and $W'$ belong to the two different Witt towers of skew-Hermitian spaces (with $\dim W \mod 2$ fixed), and $\dim W + \dim W' = 2 \dim V$. 
Then for any $\pi \in {\rm Irr}(\U(V))$,  exactly one of the two theta lifts $\Theta_{V,W,\psi}(\pi)$ and $\Theta_{V, W', \psi}(\pi)$ is nonzero.
\end{cor}

 \vskip 5pt

\subsection{\bf $\U_1 \times \U_1$} 
Let us illustrate the dichotomy theorem in the base case where $\dim V = \dim W_0 = 1$.  Let $W_0$ and $W'_0$ be the two skew-Hermitian spaces of  dimension $1$.
For any $\chi \in {\rm Irr}(\U(V)) = {\rm Irr}(E^1)$, the dichotomy theorem implies that exactly one of the theta lifts $\theta_{V, W_0, \psi}(\chi)$ or $\theta_{V, W'_0, \psi}(\chi)$ is nonzero.
Now here is an interesting question: which of these lifts is nonzero? This turns out to be a highly nontrivial and beautiful 
 result  of Moen and Rogawski \cite[Prop. 3.4]{R} and Harris-Kudla-Sweet \cite{HKS}:
  
 \begin{thm}  \label{T:epsilon}
 The theta lift $\theta_{V, W_0, \psi}(\chi)$ (with respect to splitting characters $(\chi_V, \chi_W)$) is nonzero if and only if
 \[  \epsilon(V) \cdot \epsilon(W_0) = \epsilon_E(1/2, \chi_E \cdot \chi_W^{-1} , \psi(Tr_{E/F}(\delta -))). \]
 Here, recall that $\delta$ is a  nonzero trace zero element of $E$ and the definition of the sign $\epsilon(W_0)$ depends on $\delta$. 
Moreover, $\chi_E$ is the character of $E^{\times}/F^{\times}$ defined by $\chi_E(x) = \chi(x/x^c)$ and  the local epsilon factor on the right is that defined in  Tate's thesis.
 \end{thm}

This theorem shows that the question of nonvanishing of theta lifting has deep arithmetic connections.
One way to prove this theorem is via the doubling seesaw argument, which relates the theta lifting to the doubling zeta integral of Piatetski-Shapiro and Rallis, a theory that produces  the standard L-function and epsilon factor. This is the approach of \cite{HKS}. We have no time to go into this here, but would like to mention that this doubling zeta integral plays an important role in Ellen Eischen's lectures.
\vskip 5pt
  
  After addressing nonvanishing, the next issue is that of identity: what is $\theta_{V, W_0, \psi}(\chi)$  if it is nonzero?
   Suppose we pick $\chi_V = \chi_W$ (as is allowed here). Then the splitting over the two $\U_1$'s agree (on identifying them with $E^1$) and so  the theta lifting is the identity map $\Theta(\chi) = \chi$ on its domain (i.e. outside its zero locus). With our knowledge of how the Weil representation changes when we change $(\chi_V, \chi_W)$, this allows one to figure out the general case:
\[   \Theta_{\chi_V, \chi_W, \psi}(\chi) = \chi \cdot (\chi_{W}^{-1} \chi_V \circ i) \quad \text{on its domain.} \]

\vskip 5pt

 \subsection{\bf Application to Howe-PS setting}
We shall now specialize to the particular case we are interested in. Set  $V = \langle 1 \rangle$ to be the 1-dimensional Hermitian space with form $(x,y) \mapsto xy^c$ 
so that
\[  \U(V) = E^1 \cong E^{\times}/ F^{\times}. \]
 We  consider the theta correspondence for $\U(V)$ with the two odd-dimensional Witt towers $(W_r)$ and $(W'_r)$, with $\dim W_r = \dim W'_r =2r+1$.  Because $\U(V)$ is compact, one in fact has a direct sum  decomposition:
 \[ \Omega = \bigoplus_{\chi }  \chi \otimes \Theta(\chi) \]
as $\chi$ runs over the characters of $E^1$.  Now let us record  some consequences of the general results discussed above:
\vskip 5pt

\begin{itemize}
\item For any $\chi$, $\Theta_{V,W_r, \psi}(\chi)$ and $\Theta_{V, W', \psi}(\chi)$ are irreducible or $0$. This is because, with $\U(V)$ being compact, any $\chi$ is supercuspidal.
\vskip 5pt

\item For any $\chi$ and $r >0$, $\theta_{V, W_r, \psi}(\chi)$ and $\theta_{V, W'_r,\psi}(\chi)$ are both  nonzero; this is because we are already in the stable range when $r > 0$. 
\vskip 5pt

\item  What if $r = 0$ (i.e. the dual pair $\U_1 \times \U_1$)?  This is the situation addressed by the dichotomy theorem (Corollary \ref{C:dichotomy}): exactly one of 
$\theta_{V, W_0, \psi}(\chi)$ and $\theta_{V, W'_0,\psi}(\chi)$ is nonzero.  Exactly which one is nonzero is highly non-obvious but is given by Theorem \ref{T:epsilon}. 
\vskip 5pt

\item Suppose without loss of generality that $\theta_{V, W'_0, \psi}(\chi) = 0$. Then $\Theta_{V, W'_1, \psi}(\chi)$ is supercuspidal.
\end{itemize}

 I will leave it as an exercise for the reader to deduce the above assertions from the results discussed above.
Instead, I will describe the proof of some of those results in the special case of $\U_1 \times \U_3$. This is where we do the ``dirty work", which will be formulated as a series of guided exercises below.
\vskip 10pt

\subsection{\bf Guided exercise}  \label{SS:guided-ex}
The Weil representation $\Omega$ for $\U(V) \times  \U(W_1) = \U_1 \times \U_3$ can be given a realization as follows. Let $\omega$ be the Heisenberg-Weil representation for 
\[  (\U(V) \times \U(W_0)) \ltimes H(V\otimes_E W_0) = (\U_1 \times \U_1) \ltimes H(V\otimes_E W_0)\]
where we recall that $H(V\otimes_E W_0)$ is the Heisenberg group associated to the 2-dimensional symplectic space $V \otimes_E W_0$ over $F$.
It is in fact not easy to give a concrete model for this representation (even starting with a Schrodinger model for $\Mp_2$). Then the Weil representation $\Omega$ for $\U(V) \times \U(W_1)$ can be realized on 
\[  S(E e^* \otimes  V) \otimes \omega \]
 which can be thought of as the space of Schwartz functions on the 1-dimensional $E$-vector space $E e^* \otimes V$ valued in $\omega$. One can give explicit formulas for the action of $\U(V) \times B$, where $B = TU$ is the Borel subgroup stabilizing the isotropic line $E \cdot e$ as follows.
\vskip 5pt

\begin{itemize}
\item[(a)]  for $h \in \U(V) = E^1$, 
\[  (h \cdot f) (x) = \omega(h) \left( f(h^{-1} \cdot x) \right),  \quad \text{with $x \in E$ and $f \in \Omega$.} \]
\vskip 5pt

\item[(b)]  for an element
\[  t(a, b)  = \left( \begin{array}{ccc}
 a & & \\
  & b &  \\
   &  & (a^c)^{-1}  \end{array} \right) \in T  \quad \text{with $a \in E^{\times}$ and $b \in E^1$,}  \]
one has:
   \[ (t(a,b) f)(x)  = \chi_V(a) \cdot  |a|^{1/2} \cdot  \omega(b) \left(f  (a^c \cdot x) \right), \]
   where we regard $b \in E^1$ as an element of $\U(W_0)$.  
\vskip 5pt

\item[(c)]  for an element
\[  u(0, z) = \left( \begin{array}{ccc}
1 & 0& z \\
& 1 & 0 \\
&  &  1  \end{array}  \right) \in U,  \quad \text{with $z \in F$,}\]
one has:
\[  (u(0,z) f)(x)  =  \psi(z N(x)) \cdot f(x). \]
\vskip 5pt
\item[(d)]  for an element
\[  u(y,0) =  \left( \begin{array}{ccc}
1 & * &  *  \\
& 1 &  y  \\
&  &  1  \end{array}  \right), \]
one has:
\[  (u(y,0) f) (x) =  \omega(h(xy , 0))(f(x)), \]
 where $h(xy) = (xy,0) \in H(V \otimes_E W_0)$ is regarded as an element in the Heisenberg group.
 (The two asterisks in $u(y,0)$ are determined by $y$; work out what they should be).
\end{itemize}
\vskip 10pt

Now recall that we have:
\[  \Omega =  \bigoplus_{\chi \in {\rm Irr}(E^1)}  \chi \otimes \Theta(\chi), \] 
and our goal is now to understand $\Theta(\chi)$ as much as possible. 
\vskip 5pt

Given the above information, here is the guided exercise:
\vskip 5pt
\noindent{\bf \underline{Exercise}:}
\vskip 5pt
\begin{itemize}
\item[(i)]  Let $Z = \{ u(0,z): z \in F \} \subset U$. This is the center of $U$. Compute the coinvariant space 
\[  \Omega_Z = \Omega / \langle z \cdot f - f : z \in Z, f \in \Omega \rangle, \]
as a representation of $\U(V) \times B / Z  = E^1 \times T \cdot U/Z$. Indeed, show that the natural projection $\Omega \rightarrow \Omega_Z$ is given by the evaluation-at-0 map
\[  ev_0:  \mathcal{S}(E) \otimes \omega \longrightarrow \omega. \]


 \vskip 5pt

\item[(ii)]  From the answer in (i), deduce the following:
\vskip 5pt

\begin{itemize}
\item [(a)]  For any $\chi \in {\rm Irr}(\U(V))$, $\Theta(\chi)_Z = \Theta(\chi)_U$. 
\vskip 5pt

\item[(b)]  Suppose that  $\chi$  has nonzero theta lift  $\theta_0(\chi)$ to $\U(W_0)$ with respect to $\omega$. Then $\Theta(\chi)$ is nonzero and non-supercuspidal. Indeed, one has a nonzero $T = E^{\times} \times \U(W_0)$-equivariant map 
\[  \Theta(\chi)_N \rightarrow  \chi_V|-|^{1/2}  \otimes \theta_0(\chi),   \]
so that by Frobenius reciprocity,  there is a nonzero equivariant map
\[  \Theta(\chi)  \rightarrow  {\rm Ind}_B^{\U(W)} \left( \chi_V | - |^{-1/2} \otimes \theta_0(\chi) \right) \quad \text{(normalized induction)} \]
taking note that $\delta_B^{1/2}(t(a,b)) = |a|_E$.
Hence we see that $\Theta(\chi)$ contains a constituent of the latter principal series representation, which is nontempered (since $|-|^{1/2}$ is not unitary).
\vskip 5pt

\item[(c)]   Suppose that $\chi$ has zero theta lift to $\U(W_0)$. From (i), deduce that $\Theta(\chi)$ is supercuspidal (i.e. $\Theta(\chi)_U = 0$).
\end{itemize}
\vskip 5pt

 \item[(iii)] Now compute the twisted  coinvariant space
\[  \Omega_{Z, \psi} =   \Omega / \langle z \cdot f -  \psi(z) \cdot f : z \in Z, f \in \Omega \rangle \]
as a representation of $\U(V) \times T_{\psi}$, where 
\[ T_{\psi} = \{ t(a,b):  a,b \in E^1 \}  \subset T \]
is the stabilizer of $(Z, \psi)$ in $T$. 
\vskip 5pt

\item[(iv)] Deduce from the answer in (iii) that for any $\chi \in {\rm Irr}(\U(V))$, $\Theta(\chi) \ne 0$. 
\end{itemize}
\vskip 10pt

What one sees from this guided exercise is that to understand the theta lifts $\Theta(\pi)$ (for example to detect its nonvanishing or supercuspidality), it is useful to consider various twisted coinvairant spaces $\Omega_{N, \psi}$ where $N \subset \U(W)$ is an abelian subgroup and $\psi$ is a (possibly trivial) character of $N$. 
Such twisted coinvariant spaces (or twisted Jacquet modules) are local analogs of the Fourier coefficients of modular forms. They
 are readily computable from the concrete model of the Weil representation analogous to the one above.
It also shows that the Weil representations have an inductive structure with respect to the Rallis tower. In this guided exercise, we see that we are basically reduced to the following two problems:
\vskip 5pt

\begin{itemize}
\item the irreducibility of $\Theta(\chi)$ (as given by the Howe duality theorem);
\item the understanding of the $\U_1 \times \U_1$ theta correspondence (which we  discussed earlier).
 \end{itemize}
\vskip 10pt

 \subsection{\bf Split case}
Note that the case when $E = F \times F$ is also necessary for global applications. In this case, the dual pair is ${\rm GL}_m \times {\rm GL}_n$. The Weil representation is, up to twists by 1-dimensional characters, the natural action of $\GL_m (F) \times \GL_n(F)$ on the space $\mathcal{S}(M_{m \times n}(F))$ of Schwartz functions on the space of $m \times n$ matrices.  The study of this local theta correspondence  is essentially completed in the paper \cite{Mi} of A. Minguez.   Hence we shall say no more about this case in this paper.  
\vskip 5pt

\subsection{\bf Remarks} 
We have given a discussion of the theory of classical theta correspondence which is based on reductive dual pairs in the symplectic group. But this idea is clearly more robust. One may ask:
\vskip 5pt

\begin{itemize}
\item Can one classify all reductive dual pairs $G \times H$ in any simple Lie group $E$, as opposed to just for $E = \Sp_{2n}$?

\item Is there an understanding of the smallest infinite-dimensional representation $\Omega$ of any such $E$?

\item If so, when one pulls back $\Omega$ to $G \times H$, does one obtain a  transfer or lifting of representations analogous to those described in this lecture? 
\end{itemize}

These questions started to be explored in the mid-1980's. Reductive dual pairs have been classified on the level of Lie algebras by Rubenthaler. The construction and classification of the so-called minimal representations of a simple Lie group $E$ was begun by Kostant, Vogan, Kazhdan, Savin, Torasso and others; see \cite{GS}. Finally the study of the resulting theta correspondence began in the 1990's but the  theory is not as systematic as the classical case. It is only recently that one has somewhat complete results in several families of examples. In the project for this course, you will work with a particular instance of this exceptional theta correspondence.
 
\newpage

\section{\bf Lecture 3: Global Theta Correspondence}
In this third lecture, we will discuss the global theta correspondence. For the case of unitary groups, a nice reference is the article \cite{Ge} of Gelbart. 
We will see that almost all of the considerations and constructions of the previous lecture make sense in the global setting, once they are appropriately construed. We will  work over a number field $k$   with associated local field $k_v$ for each place $v$ of $k$ and with adele ring $\A = \prod_v' k_v$. We fix a quadratic field extension $E/k$ and consider a pair of a Hermitian space $V$ and a skew-Hermitian space $W$ relative to $E/k$.
\vskip 5pt

\subsection{\bf Basic idea.}
Let us return to the basic idea of Lecture 2: in the local setting, the Weil representation allows one to define the local theta lifting
\[  \theta : {\rm Irr}(\U(V)) \longrightarrow {\rm Irr}(\U(W)) \cup \{0 \}. \]
If $\U(V) \times \U(W)$ is in the stable range (with $V$ smaller), one even has
\[  \theta : {\rm Irr}(\U(V)) \longrightarrow {\rm Irr}(\U(W)). \]
 In the global setting, one might imagine that by taking (restricted) tensor product taken over all places $v$ of a number field $k$, one gets
\[  \theta :  {\rm Irr}(\U(V)(\A)) \longrightarrow {\rm Irr}(\U(W)(\A)). \]
 This is the case, but we are interested not in the lifting of abstract irreducible representations, but rather of cuspidal automorphic representations. Cuspidal automorphic representations
are representations which are realized in the space of cuspidal automoprhic forms (i.e. functions on $[G]$). So what we need is a map
\[  \{\text{Cusp forms on $\U(V)$} \} \longrightarrow \{ \text{Automorphic  forms on $\U(W)$} \}. \]

\vskip 5pt
We are thus interested in procedures which allow one to lift functions on a (measure) space $X$ to functions on another space $Y$. A standard such procedure is via a kernel function $K$, i.e. a function $K : X \times Y \rightarrow \C$. Given such a function $K$, one gets a linear map 
\[  T_K : \mathcal{C}(X) \longrightarrow  \mathcal{C}(Y) \]
defined by
\[  T_K(f)(y) = \int_X  K(x,y) \cdot f(x) \, dx, \]
assuming convergence is not an issue.
\vskip 5pt

Let's apply this simple idea to our setting. Recall from Lecture 2 that  we have
\[  \tilde{\iota}: \U(V) \times \U(W) \longrightarrow  \Mp(V \otimes_E W). \]
We shall see that the (global) Weil representation $\Omega$ is automorphic on $\Mp(V \otimes_E W)$, i.e. there is an equivariant map
\[  \theta:  \Omega \longrightarrow \mathcal{A}(\Mp(V \otimes_E W)). \]
So for any $\phi \in\Omega$, we have an automorphic form $\theta(\phi)$ on $\Mp(V \otimes_E W)$: these are the theta functions. 
Pulling back $\theta(\phi)$ by $\tilde{\iota}$, we may regard $\theta(\phi)$  as a function on $[\U(V)] \times [\U(W)]$.
We can thus use $\theta(\phi)$ as a kernel function to transfer functions on $[\U(V)]$ to functions on $[\U(W)]$. In other words, each $\theta(\phi)$ gives a linear map
\[   \theta_{\phi}: \{\text{ Cusp forms on $\U(V)$} \} \longrightarrow \{ \text{Automoprhic forms on $\U(W)$} \}. \]
As we consider all these $\theta_{\phi}$ together, we have a map 
\[ \begin{CD}
 \{ \text{Cuspidal automorphic representations of $\U(V)$} \} \\
 @VVV  \\
 \{\text{Automorphic representations of $\U(W)$} \}. \end{CD}  \]
\vskip 5pt

\subsection{\bf Adelic metaplectic groups}
We shall now give more precise formulation of the above basic idea. Fix a non-trivial additive character $\psi = \prod_v \psi_v$ of $F \backslash \A$. 
Suppose that $W$ is a symplectic vector space over $k$. Then for each $v$, we have seen the metaplectic group
\[  \begin{CD}
1 @>>> S^1 @>>> \Mp(W_v) @>>>\Sp(W_v) @>>> 1 \end{CD} \]
For almost all $v$, it is known that the covering splits uniquely over the hyperspecial maximal compact subgroup $K_v$, so that we may regard $K_v$ as an open compact subgroup opf $\Mp(W_v)$. 
Then one can form the restricted direct product:
\[  {\prod}'_v \Mp(W_v) \quad \text{ (with respect to the family $\{K_v\}$).} \]
This contains as a central subgroup $\bigoplus_v S^1$. If we quotient out the restricted direct product above by the central subgroup
\[   Z = \{ (z_v) \in \bigoplus_v S^1:  \prod_v z_v = 1 \} \]
we get the adelic metatplectic group
\[  \begin{CD}
1@>>> S^1 @>>> \Mp(W)(\A) @>>> \Sp(W) @>>> 1. \end{CD} \]
Note that though we use the notation $\Mp(W)(\A)$, $\Mp(W)$ is not an algebraic group and we are not taking the group of adelic points of an algebraic group.  
\vskip 5pt

An important property of this adelic metaplectic cover is that it splits  (canonically) over the group $\Sp(W)(k)$ of $k$-rational points, so that one can canonically regard $\Sp(W)(k)$ as a subgroup of $\Mp(W)(\A)$. As a result one can consider the automorphic quotient
\[  [\Mp(W)] = \Sp(W)(k) \backslash \Mp(W)(\A) \]
and introduce the space of genuine automorphic forms on $\Mp(W)(\A)$: these are the automorphic functons 
\[  f: [\Mp(W)] \longrightarrow \C \]
such that for all $z \in S^1$,
\[   f(zg) = z \cdot f(g). \]
\vskip 5pt

\subsection{\bf Global Weil representations}
We may consider the global Weil representation
\[ \omega_{\psi} :=  {\bigotimes}'_v \omega_{\psi_v} \quad \text{of \, ${\prod}_v' \Mp(W_v)$.} \]
This factors to a representation $\omega_{\psi}$ of $\Mp(W)(\A)$ (for clearly $Z$ acts trivially).  
If $W = X \oplus Y$ is a Witt decomposition, we have seen that for each $v$, $\omega_{\psi_v}$ can be realized on $\mathcal{S}(Y_v)$. Hence, $\omega_{\psi}$ can be realized on 
\[  {\bigotimes}_v' \mathcal{S}(Y_v) =  \mathcal{S}(Y_{\A}). \]
In other words, $\omega_{\psi}$ is realized on a very concrete space of functions. 
\vskip 5pt

\subsection{\bf Theta functions}
It turns out that one has an $\Mp(W)(\A)$-equivariant map
\[  \theta :  \mathcal{S}(Y_{\A}) \longrightarrow  \mathcal{A}(\Mp(W)) \]
 defined by averaging over the $k$-rational points of $Y$:
 \[  \theta(f)(g) = \sum_{y \in Y_k}  (\omega_{\psi}(g) \cdot f)(y). \]
 The fact that $\theta(f)$ is left-invariant under $\Sp(W)(k)$ is a consequence of the Poisson summation formula.
The functions $\theta(f)$ are called theta functions.
 \vskip 5pt

The map $\theta$ is non-injective. More precisely, since 
\[   \omega_{\psi_v} = \omega_{\psi_v}^+ \oplus \omega_{\psi_v}^-  \quad \text{for each $v$,} \]
we see that as an abstract representation,
\[  \omega_{\psi}   = \bigoplus_S \omega_{\psi, S}  \]
where
\[  \omega_{\psi,S} = (\bigotimes_{v \in S} \omega_{\psi_v}^-) \otimes (\bigotimes_{v \notin S} \omega_{\psi_v}^+) \]
for finite subsets $S$ of places of $k$.  Then 
\[  {\rm Ker}(\theta) = \bigoplus_{\text{$|S|$ odd}} \omega_{\psi, S}, \]
and we have an injective map (still denoted by $\theta$)
\[ \theta:  \bigoplus_{\text{$|S|$ even}} \omega_{\psi,S} \hookrightarrow \mathcal{A}(\Mp(W)) \]
upon restriction to the ``even" subspace above.

 \vskip 5pt

\subsection{\bf Pulling back} 
Now  let us return to our reductive dual pair $\U(V) \times \U(W)$ over $k$.
Recall that we have
\[  \iota: \U(V)(\A) \times \U(W)(\A) \longrightarrow \Sp(V \otimes_E W)(\A). \]
If we fix a pair of automorphic  characters $(\chi_V, \chi_W)$ of $E^{\times}$ with
\[  \chi_V |_{\A^{\times}} = \omega_{E/k}^{\dim_E V} \quad \text{and} \quad \chi_W|_{\A^{\times}} = \omega_{E/k}^{\dim_E W}, \]
 then as in the local case, one obtains an associated lifting
\[  \tilde{\iota}_{\chi_V, \chi_W, \psi}: \U(V)(\A) \times \U(W)(\A) \longrightarrow  \Mp(V \otimes W)(\A). \]
Using this lifting, we may pullback the global Weil representation $\omega_{\psi}$ to obtain the Weil representation
\[  \Omega = \Omega_{V,W,\chi_V, \chi_W, \psi} \quad \text{  of $\U(V)(\A) \times \U(W)(\A)$.} \]
Moreover, the lifting $\tilde{\iota}$ sends the group $\U(V)(k) \times \U(W)(k)$  into $\Sp(V \otimes W)(k) \subset \Mp(V \otimes W)(\A)$. Hence, the pullback of a function in $\mathcal{A}(\Mp(V \otimes W))$ by $\tilde{\iota}$ gives a smooth function on $[\U(V)] \times [\U(W)]$. Thus we have
\[  \theta : \Omega \longrightarrow \mathcal{A}(\Mp(V \otimes W)) \longrightarrow \mathcal{C}^{\infty} ( [\U(V) \times \U(W)]). \]  
\vskip 5pt

\subsection{\bf Global theta liftings}
Now  for $f \in \mathcal{A}_{cusp}(\U(V))$ and $\varphi \in \Omega$, we set:
\[ \theta(\varphi, f)(g) = \int_{[\U(V)]}  \theta(\varphi)(g, h) \cdot \overline{f(h)} \, dh. \] 
The cuspidality of $f$ implies that the integral above converges absolutely (because of the rapid decay of $f$ modulo center).
Then $\theta(\varphi, f)$ is an automorphic form on $\U(W)$. 

\vskip 5pt
Suppose that $\pi \subset \mathcal{A}_{cusp}(\U(V))$ is an irreducible cuspidal automorphic representation. Let
\[  \Theta(\pi) = \langle \theta(\varphi,f): f \in \pi, \, \varphi \in \Omega \rangle \subset \mathcal{A}(\U(W)).  \]
This is a $\U(W)(\A)$-submodule of $ \mathcal{A}(\U(W))$ and we call it the  {\bf global theta lift} of $\pi$.
\vskip 10pt

\subsection{\bf Questions}

 The main questions concerning global theta lifting are analogs of those in the local case:
 \vskip 5pt
 
 \begin{itemize}
 \item Is $\Theta(\pi)$ nonzero?
 
 \vskip 5pt
 
 \item Is $\Theta(\pi)$ contained in the space of cusp forms? 
\end{itemize}

In addition, we can ask:
\begin{itemize}
\item For $\pi = \otimes_v' \pi_v$, how is the global $\Theta(\pi)$ related to the local theta liftings $\theta(\pi_v)$ for all $v$? 
\end{itemize}
We shall begin by addressing this issue of local-global compatibility.
\vskip 5pt

\subsection{\bf Compatibility with local theta lifts}
How is the representation $\Theta(\pi)$ related to the abstract irreducible representation $\Theta^{abs}(\pi) : = \otimes_v \theta(\pi_v)$?
We have:

\begin{prop}
 Suppose that $\Theta(\pi)$ is non-zero and is contained in the space $\mathcal{A}_2(\U(W))$ of square-integrable automorphic forms on $\U(W)$. 
Then $\Theta(\pi) \cong \Theta^{abs}(\pi)$.
\end{prop}

\proof
Since  $\Theta(\pi) \subset \mathcal{A}_2(\U(W))$, it is semisimple. Let $\sigma$ be an irreducible summand of $\Theta(\pi)$. Then consider the linear map 
\[  \Omega   \otimes \pi^{\vee} \otimes \sigma^{\vee} \longrightarrow \mathbb{C} \] 
defined by: 
\[ \varphi \otimes \overline{f_1} \otimes
\overline{f_2} \mapsto \int_{[\U(W)]} \theta(\varphi, f_1)(g) \cdot \overline{f_2(g)} \, dg. \] 
This map is non-zero and
$\U(V) \times \U(W)$-invariant. Thus it gives rise to a non-zero equivariant map 
\[ \Omega  \longrightarrow \pi \otimes \sigma, \] 
and thus for all $v$, a non-zero $\U(V_v) \times \U(W_v)$-equivariant map 
\[ \Omega_v \longrightarrow \pi_v \otimes \sigma_v. \] 
In other words, we must have 
\[ \sigma_v \cong \theta(\pi_v). \] 
Hence, $\Theta(\pi)$ must be an isotypic sum of $\Theta^{abs}(\pi)$.  
Moreover, the multiplicity-one statement in the Howe duality theorem implies that 
\[ \dim \Hom_{\U(V)(\A) \times \U(W)(\A)} (\Omega, \pi \otimes \Theta^{abs}(\pi)) = 1. \] 
Thus $\Theta(\pi)$ is in fact irreducible and isomorphic to $\Theta^{abs}(\pi)$.
\qed

\vskip 5pt

\subsection{\bf  Cuspidality and Nonvanishing}
As in the local case, it is useful to consider a Rallis tower of theta lifitngs, corresponding to a Witt tower $W_r = W_0 \oplus \mathbb{H}^r$ of skew-Hermitian spaces, with $W_0$ anisotropic. 
One has the analogous results:
\vskip 5pt

\begin{prop}  \label{P:rallis}
Let $\pi$ be a cuspidal automorphic representation of $\U(V)$, and consider its global theta lift $\Theta_{V, W_r,\psi}(\pi)$ on $\U(W_r)$ (relative to a fixed pair $(\chi_V, \chi_W)$). 
Then one has:
\vskip 5pt

(i) There is a smallest  $r_0 = r_0(\pi) \leq \dim V$ such that $\Theta_{V, W_{r_0}, \psi}(\pi) \ne 0$.  Moreover,  $\Theta_{V, W_{r_0}, \psi}(\pi) $ is contained in the space of cusp forms.

\vskip 5pt

(ii) For all $r >  r_0$, $\Theta_{V, W_{r}, \psi}(\pi)$ is nonzero and noncuspidal.

\vskip 5pt

(iii) For all $r \geq  \dim V$, $\Theta_{V, W_{r_0}, \psi}(\pi)  \subset \mathcal{A}_2(\U(W))$, unless $ r = \dim V$, $r_0 = 0$ and $W_0 = 0$.  
\end{prop}
As in the local case, we call $r_0= r_0(\pi)$ the first occurrence of $\pi$ in the relevant Witt tower, and we call the range where $r \geq \dim V$ the stable range.  
\vskip 5pt

Let us give some indication of the proof of parts (ii) and (iii) of this proposition, assuming the existence of $r_0$ as in (i).  The key input is Rallis' computation of the constant term of a global theta lift \cite[Thm. I.1.1]{Ra}. More precisely,  for a maximal parabolic subgroup $P_j = M_j \cdot N_j$ of $\U(W_r)$ with Levi factor $M_j \cong \GL_j(E) \times \U(W_{r-j})$ (for $ 1 \leq j \leq r$), one considers the normalized constant term  
\[  R_{P_j}( \Theta_{V,W_r, \psi}(\pi)) := \delta_{P_j}^{-1/2}  \cdot \langle  \theta(\varphi,f)_{N_j} : f \in \pi, \,  \varphi \in \Omega_{V, W_r,\psi} \rangle \]
of $\Theta_{V,W_r, \psi}(\pi)$ along $N_j$: this is a space of automorphic forms on $M_j$ and is the global analog of the normalized Jacquet module in the local setting.  Rallis' result \cite[Thm. I.1.1]{Ra} states that  
\[  R_{P_j}( \Theta_{V,W_r, \psi}(\pi)) = |\det|^{\frac{\dim V - \dim W_r + j}{2}}  \boxtimes    \Theta_{V, W_{r-j}, \psi}(\pi) \]
where $\det$ refers to the determinant character of $\GL_j(\A_E)$. Since $\Theta_{V, W_{r-j}, \psi}(\pi) = 0$ if $r-j < r_0$  (by the definition of $r_0$), we deduce:
\vskip 5pt

\begin{itemize}
\item All constant terms of $\Theta_{V, W_{r_0}, \psi}(\pi) $ vanish, so that $\Theta_{V, W_{r_0}, \psi}(\pi)$ is  cuspidal (and nonzero).
\vskip 5pt

\item For $r > r_0$, $\Theta_{V, W_r, \psi}(\pi) $ is noncuspidal (and hence nonzero), since its normalized constant term with respect to $P_{r-r_0}$ is nonzero.
Indeed,  its normalized constant term along $P_{r-r_0}$ is
\[  |\det|^{\frac{\dim V  - \dim W_0 - (r+r_0)}{2}}  \boxtimes   \Theta_{V, W_{r_0}, \psi}(\pi). \]
\vskip 5pt

\item  It follows that for $r > r_0$,  the cuspidal support of $\Theta_{V, W_r, \psi}(\pi)$ is the parabolic subgroup $Q_{r_0}$ with Levi factor $\GL_1(E)^{r-r_0} \times \U(W_{r_0})$. 
Further,   $\Theta_{V, W_r, \psi}(\pi)$ has a unique cuspidal exponent with respect to $Q_{r_0}$, given by the following  character of $\GL_1(E)^{r-r_0}$:
\[  | -|^{\frac{\dim V  - \dim W_0 - (r+r_0)}{2}} \cdot \left( |-|^{-\frac{j-1}{2}} \times |-|^{-\frac{j-3}{2}} \times .....\times |-|^{\frac{j-1}{2}} \right). \]
When $r \geq  \dim V$, the exponent 
\[ \dim V  - \dim W_0 - (r+r_0)   < 0,  \]
 unless $r = \dim V$, $r_0 = 0$ and $W_0 =0$. 
 Since this exponent is negative, it follows from the square-integrability criterion of Jacquet \cite[Lemma I.4.11]{MW} that $\Theta_{V, W_r, \psi}(\pi)$ is square-integrable.
 
\end{itemize}
\vskip 10pt

\subsection{\bf The case of $\U_1 \times \U_1$} 
In parallel with the local setting, we may consider the theta lift for the basic case of $\U(V) \times \U(W)$ with $\dim V = \dim W = 1$. We have seen that the nonvanishing of local theta lifts is controlled by a local root number. Here is the global theorem \cite{Y}:
\vskip 5pt

\begin{thm}
Let $\chi$ be an automorphic character of $\U(V) = E^1$. Its global theta lift $\Theta_{\chi_V, \chi_W, \psi}(\pi)$ on $\U(W)$ is nonzero if and only if :
\vskip 5pt

\begin{itemize}
\item [(i)] for each place $v$, the local theta lift $\theta_{\chi_{V,v}, \chi_{W,v}, \psi_v}(\chi)$ is nonzero;
\item[(ii)]  the global L-value 
\[  L(1/2, \chi_E \chi_W^{-1}) \ne 0. \]
\end{itemize}
\end{thm}
Note that under condition (i), our local theorem for $\U_1 \times \U_1$  implies that
\[  \epsilon(1/2, \chi_{E,v} \chi_{W,v}^{-1}, \psi_v({\rm Tr} (\delta-)))   = \epsilon(V_v) \cdot \epsilon(W_v)  \quad \text{for all $v$,} \]
where $\delta$ is a trace $0$ element of $E$.  On taking product over all places $v$, we see that
\[  \epsilon(1/2, \chi_E \chi_W^{-1})  = 1  \]
since $\prod_v \epsilon(V_v) = 1 = \prod_v \epsilon(W_v)$. Hence, there is a chance that the global L-value in (ii) is nonzero!
\vskip 10pt

The proof of this Theorem is a global analog of the local theorem in the $\U_1 \times \U_1$ case, via the doubling seesaw and doubling zeta integral. One then gets the  Rallis inner product formula, which relates the Petersson inner product of two global theta lifts and the special L-value. This result thus gives an interpretation for the nonvanishing of the central L-value.
 
\vskip 10pt
\subsection{\bf The case of $\U_1 \times \U_3$}
Let us now specialize to the case of interest, where $\dim V = 1$ and $\dim W  = \dim W_1= 3$. 
Let $\chi$ be an automorphic character of $E^1$. Then $\chi$ is cuspidal since $E^1$ is anisotropic. We can thus apply the above general results to conclude:
 
\begin{cor} \label{C:global}
$\Theta_{V, W_1}(\chi)$ is nonzero and square-integrable. It is cuspidal if and only if the global theta lift of $\chi$ to $\U(W_0) = \U_1$ is zero. 
\end{cor}

\vskip 10pt 
\subsection{\bf Counterexample to Ramanujan}
Using all the results we have seen, we can now construct a counterexample to the naive Ramanujan conjecture:
\vskip 5pt

\begin{itemize}
\item Consider the global theta correspondence for $\U(V) \times \U(W_1)$ with $\dim V  = 1$ and $\dim W_1  =3$, and take $\chi_V = \chi_W$ and  $\chi = 1$ (the trivial character of $\U(V)(\A)$).
The global theta lift $\Theta_{V,W,\psi}(1)$ is a nonzero irreducible square-integrable automorphic representation. 
\vskip 5pt

\item If $1$ has zero global theta lift   to the lower step   $\U(W_0) = \U_1$ of the Rallis tower, then
$\Theta(1)$ is cuspidal irreducible. Moreover,  by  part (ii)(b) of  the guided exercise in \S \ref{SS:guided-ex}, one sees that for almost all $v$, $\Theta(1)_v \cong \theta(1_v) $ is an unramified representation belonging to the principal series 
\[  {\rm Ind}_{B_v}^{\U(W_v)}\left(\chi_{V,v} |-|_v^{-1/2} \otimes  1_v\right) \]
and thus is nontempered. Hence $\Theta(1)$ would be a counterexample  to the naive Ramanujan conjecture.
\vskip 5pt

\item  On the other hand, if $1$ has nonzero theta lift to $\U(W_0)$,  then we may select two places $v_1$ and $v_2$ and replace each of $W_{v_1}$ and $W_{v_2}$ by the other local skew-Hermitian space. In other words, we can find a global skew-Hermitian space $W'$ such that
\[  W'_{v_1} \ne W_{v_1}  \quad \text{and} \quad  W'_{v_2} \ne W_{v_2} \]
but
\[  W'_v \cong W_v  \quad \text{for all $v \ne v_1, v_2$.}  \]
Such a global skew-Hermitian space exists,  by our classification of global Hermitian spaces.
\vskip 5pt

Now consider the global theta lift $\Theta_{V, W'}(1)$ on $\U(W')$ and observe that 
\[ \theta_{V_{v_1}, W'_{v_1}}(1) = 0 = \theta_{V_{v_2}, W'_{v_2}}(1)  \]
because of dichotomy. Hence we can repeat the above argument in replacing $W$ by $W'$.
\end{itemize}
 \vskip 5pt
 
 \subsection{\bf Guided exercise} 
 As in the local case, it will be instructive to carry out some of the computations in proving some of the above results, at least in our setting of $\U_1 \times \U_3$. 
 The guided exercise below is the global analog of the local guided exercise in Lecture 2, following the various notation there.
 
 \vskip 5pt
 
 To do the exercise,  it is necessary to write down the theta function $\theta(\phi)$ (for $\phi \in \Omega$) as explicitly as possible. Recall that $\Omega$ is realized on
\[  \mathcal{S}(\A_E e^* \otimes v_0) \otimes \omega \]
where $\omega$ is the global Heisenberg-Weil representation of 
\[  (\U(V) \times \U(W_0)) \ltimes H(V \otimes_E W_0). \]
 We have an automorphic functional
\[  \theta_0:  \omega \longrightarrow \mathcal{A}(\Mp_2)  \rightarrow \C \]
where the last arrow is the evaluation at $1$.  The automorphic realization
\[  \theta: \Omega = \mathcal{S}(\A_E) \otimes \omega  \longrightarrow \mathcal{A}(\U(V) \times \U(W)) \]
is then given by
\[  \theta ( f) (g)  = \sum_{x \in E}  \theta_0 \left( \Omega(g)(f)(x) \right). \]

 Now for the exercise:
 \vskip 5pt
 
 \noindent{\bf \underline{Exercise}:}
 \begin{itemize}
 \item For a cusp form $f$ on $\U(V)$ and $\phi \in \Omega$, compute:
 \[  \theta(\phi, f)_Z(g) := \int_{[Z]} \theta(\phi,f)(zg) \, dz  \]
 and
  \[  \theta(\phi,f)_U(g) = \int_{[U]} \theta(\phi,f)(ug) \, du, \]
 where $Z = \{ u(0,z): z \in k\}$ is the center of the unipotent radical $U$ of $B$.

 \vskip 5pt
 
 \item Likewise, compute:
 \[  \theta(\phi, f)_{Z,\psi}(g) := \int_{[Z]} \overline{\psi(z)} \cdot \theta(\phi,f)(zg) \, dz  \]

\item From these computations, deduce as much of Cor. \ref{C:global} as possible.
\end{itemize}
\vskip 5pt
 \newpage

\section{\bf Lecture 4: Arthur's Conjecture}
In this final lecture, we will discuss an influential conjecture of Arthur \cite{Ar1} which explains and classifies all possible failures of the naive Ramanujan conjecture. We will then illustrate with 
some examples, including the Howe-PS example discussed earlier, the Saito-Kurokawa example for $\PGSp_4$ and the case of the exceptional group $G_2$.  
\vskip 5pt

\subsection{\bf A basic hypothesis}
In the formulation of Arthur's conjecture, one needs to make a (serious) assumption: 
\vskip 5pt

\noindent{\bf (Basic Hypothesis):} There is a topological group $L_k$ (depending only on the number field $k$) satisfying the following 
properties:
\vskip 5pt

\begin{itemize}
\item the identity component $L_k^0$ of $L_k$ is compact and the group of components $L_k / L^0_k$ is isomorphic to the Weil 
group $W_k$ of $k$;
\vskip 5pt

\item for each place $v$, there is a natural conjugacy class of embeddings $L_{k_v} \hookrightarrow L_k$, where $L_{k_v}$ is the Weil 
group if $k_v$ is archimedean, and the Weil-Deligne group $W_{k_v} \times \SU_2(\C)$ if $k_v$ is non-archimedean.
\vskip 5pt

\item there is a natural bijection between the set of isomorphism classes of irreducible representations of $L_k$ of dimension $n$
and the set of cuspidal representations of $GL_n(\A)$. 
\end{itemize}

\vskip 5pt
This assumption is basically the main conjecture in the Langlands program for $GL_n$.
We can view it as a classification of the cuspidal representaitons of $\GL_n$ in terms of irreducible $n$-dimensional Galois representations.
\vskip 5pt

\subsection{\bf Arthur's conjecuture}
Arthur's conjecture is a classification of the constituents of $\mathcal{A}_2(G)$, i.e. a classification of the square-integrable automorphic representations of $G$.
This classification proceeds in two steps. The first step is approximately the classification of these constituents up to  near equivalence (this is not entirely true, but for the groups discussed here, it is expected to be so). Here, we say that two representations $\pi_1 = \otimes'_v \pi_{1,v}$ and $\pi_2 = \otimes'_v \pi_{2,v}$ of $G(\A)$ are {\em nearly equivalent} if for almost all places 
$v$, $\pi_{1,v}$ and $\pi_{2,v}$ are isomorphic: this is an equivalence relation.
\vskip 5pt

More precisely, Arthur speculated that there is a decomposition
\[  \mathcal{A}_2(G)= \bigoplus_{\psi} \mathcal{A}_{2,\psi}, \]
where each $\mathcal{A}_{2,\psi}$ is a near equivalence class, and the direct sum runs over 
equivalence classes of discrete A-parameters $\psi$. 
For a split group $G$, an A-parameter is an  admissible map
\[  \psi:L_k\times SL_2(\C) \longrightarrow G^{\vee} \]
where  $G^{\vee}$ is the complex Langlands dual group of $G$. One property of being admissible is that $\psi(L_k)$ should be bounded in $G^{\vee}$. 
An A-parameter is discrete if the centralizer group 
\[  S_{\psi} := Z_{G^{\vee}}(\psi) / Z(G^{\vee}) \]
 is finite.
 
 \vskip 5pt

The second step of the classification is to describe the decomposition of each $\mathcal{A}_{2,\psi}$.  This has a local-global structure. That is, for each place $v$, there will be a finite set of unitary representations of $G(k_v)$ (depending on the A-parameter $\psi$). If we pick an element $\pi_v$ from each of these finite sets at each place $v$, we may form the (restricted) tensor product $\pi := \otimes_v \pi_v$, which is a representation of $G(\A)$. Then $\mathcal{A}_{2,\psi}$ is the sum of such representations, with appropriate multiplicities. 
Let us now be more precise.
\vskip 5pt

\subsection{\bf Local A-packets.} The global A-parameter $\psi$  gives rise (by restriction) to a local A-parameter 
\[ \psi_v: L_{k_v}\times SL_2(\C)\longrightarrow G^{\vee}\]
for each place $v$ of $k$. 
Set 
\[  S_{\psi_v} = \pi_0 \left(Z_{G^{\vee}}(\psi_v)/ Z(G^{\vee})\right). \]
This is the local component group of $\psi_v$.
To each irreducible representation $\eta_v$ of
$S_{\psi_v}$, Arthur speculated that one can attach
a unitarizable finite length  (possibly reducible,
possibly zero) representation $\pi_{\eta_v}$ of $G(k_v)$. The set
$$A_{\psi_v}=\left\{ \pi_{\eta_v}\, :\, \eta_v \in {\rm Irr}(S_{\psi_v})
\right\} $$
is called a local A-packet. Among other things, it is required that 
\vskip 5pt

\begin{itemize}
\item for almost all $v$, $\pi_{\eta_v}$ is irreducible and unramified if $\eta_v$
is the trivial character $1_v$. In fact, for almost all $v$, $\pi_{\eta_v}$ is
the unramified representation whose Satake parameter is:
\[ s_{\psi_v} = \psi_v \left( {\rm Frob}_v \times \left( \begin{array}{cc}
q_v^{1/2} & \\ 
  &  q_v^{-1/2} \end{array} \right) \right), \]
where ${\rm Frob}_v$ is a Frobenius element at $v$ and $q_v$ is the number of elements of the 
residue field at $v$.
\end{itemize}
\vskip 5pt

\subsection{\bf Global A-packets.}
With the local packets $A_{\psi_v}$ at hand, we may define the global A-packet by:
\[  A_\psi= \left \{ \pi =\otimes_v \pi_{\eta_v} :
 \pi_{\eta_v} \in A_{\psi_v} \,  \eta_v=1_v
 \,\ \text{for almost all $v$}
\right\}.\]
It is a set of nearly equivalent representations of $G(\mathbb{A})$, indexed by the 
irreducible representations of the compact group 
\[  \mathcal{S}_{\psi,\mathbb{A}} := \prod_v S_{\psi,v}. \]
 If
$\eta = \bigotimes_v \eta_v$ is an irreducible character of 
$\mathcal{S}_{\psi,\mathbb{A}}$, then we may set
\[ \pi_{\eta} = \bigotimes_v \pi_{\eta_v}. \]
This is possible because for almost all $v$, 
$\eta_v = 1_v$ and $\pi_{1_v}$ is required to be 
unramified by the above. 
\vskip 10pt

\subsection{\bf Multiplicity formula.}
The space $\mathcal{A}_{2, \psi}$ will be the sum of the 
elements of $A_{\psi}$ with some multiplicities. 
More precisely, note that there is a natural map
\[  S_{\psi} \longrightarrow S_{\psi,\A}. \]
Arthur attached to $\psi$ a quadratic character $\epsilon_{\psi}$ of 
$S_{\psi}$ (whose definition is given below). Now if $\eta$ is an irreducible character of $S_{\psi,\A}$, we set
\[ m_{\eta} = \frac{1}{\# S_{\psi}} \cdot \left( \sum_{s \in S_{\psi}} \epsilon_{\psi}(s) 
\cdot \eta(s) \right).\]
Then Arthur conjectures that  
\[  \mathcal{A}_{\psi} \cong  \bigoplus_{\eta} m_{\eta}\pi_{\eta}. \]
\vskip 10pt

\subsection{\bf The character $\epsilon_{\psi}$}
The definition of the quadratic character $\epsilon_{\psi}$ is quite subtle. For a discrete $\psi$, one considers the adjoint action of 
\[  S_{\psi} \times L_k \times \SL_2(\C)  \quad \text{on ${\rm Lie}(G^{\vee})$ via $\psi$} \] 
and decomposes it into the direct sum of irreducible summands, each of which has the form
\[  \eta \otimes \rho \otimes S_r \]
where $S_r$ denotes the r-dimensional irreducible representation of $\SL_2(\C)$. 
Note that this is an orthogonal representation, since the adjoint representation has a nondegenerate invariant  symmetric bilinear form (e.g. the Killing form if $G$ is semisimple).
\vskip 5pt

We consider only those irreducible components $\eta \otimes \rho \otimes S_k$ satisfying the following properties:
\vskip 5pt

\begin{itemize}
\item $\eta \otimes \rho \otimes S_k$ is orthogonal;
\item $k$ is even, so that $S_k$ is a symplectic representation.
\item $\rho$ is symplectic, and 
\[  \epsilon(1/2, \rho)  = -1. \]
\end{itemize}
The conditions above imply that $\eta$ is orthogonal, so that $\det(\eta)$ is a quadratic character of $S_{\psi}$. Let $\mathcal{T}$ be the set of irreducible summands satisfying these conditions.  Then we set
\[  \epsilon_{\psi}(s) = \prod_{\tau \in \mathcal{T}}  \det(\eta)(s)   \quad \text{for $s \in S_{\psi}$.} \]
\vskip 5pt

As an example, suppose that $\psi$ is trivial on $\SL_2(\C)$. Then the set $\mathcal{T}$ is empty (since the only $S_k$ which occurs above is the trivial representation $S_1$).  

\vskip 10pt

\subsection{\bf Tempered and non-tempered parameters.}
An A-parameter $\psi$ is called {\sl tempered} if $\psi$ is 
trivial when restricted to
$SL_2(\C)$. In this case, the representations in $A_\psi$ are conjectured
to be tempered (this corresponds to the boundedness condition on $\psi(L_k)$). 
A non-tempered $A$-parameter  is one for which $\psi(\SL_2(\C))$ is not trivial.
 In this case, for almost all $v$, the unramified representation $\pi_{1_v}$ (which has Satake parameter $s_{\psi_v}$) is nontempered. 
 \vskip 5pt
 
 Thus, according to Arthur's conjecture, the cuspidal representations in
 \[  \bigoplus_{\text{nontempered $\psi$}} \mathcal{A}_{\psi} \]
 are precisely those which violate the naive Ramanujan conjecture. 
 On the other hand, the cuspidal representations in $\mathcal{A}_{\psi}$ for tempered $\psi$ are all tempered and the representation $\pi_1 =\otimes_v \pi_{1_v}$ should be globally generic. In this sense, Arthur's conjecture provides an explanation and classification of nontempered cusp forms. 
 \vskip 5pt
 
 Though the group $L_k$ is conjectural, we really only need its irreducible representations in formulating Arthur's conjecture for classical groups. As such, under our (basic hypothesis), we can replace all occurrences of ``an irreducible $n$-dimensional representation of $L_k$ by ``an irreducible cuspidal representation of $\GL_n$".  Then we may view Arthur's conjecture as a description of $\mathcal{A}_2(G)$ in terms of cuspidal representations of $\GL$'s.   Understood in this way, when  $G$ is a quasi-split classical group, Arhtur's conjecture has been verified in the works of Arthur \cite{Ar2} and Mok \cite{M}.

 \vskip 5pt
 
 \subsection{\bf Connection with theta correspondence}
 Here is a natural question one can ask concerning theta correspondence and Arthur's conjecture. 
 We have seen in Lecture 3 that when $\U(V) \times \U(W)$ is a dual pair in the stable range (with $V$ the smaller space, and barring an unfortunate boundary case), then for any cuspidal representation $\pi$ of $\U(V)$, its global theta lift $\Theta(\pi)$ on $\U(W)$ is a nonzero  irreducible summand of $\mathcal{A}_2(\U(W))$. Since square-integrable automorphic representations have A-parameters, it is natural to ask how the A-parameters of $\Theta(\pi)$ and $\pi$ are related. 
 \vskip 5pt
 
  If we view A-parameters as representing near equivalence classes, answering this question is  about identifying the local theta lifts of unramified representations and then detecting the automorphy of the family of unramified local theta lifts. This line of reasoning leads to:
 
 \vskip 5pt
 
 \begin{con}[Adam's conjecuture]
 If $\pi \subset \mathcal{A}_{cusp}(\U(V))$ has A-parameter $\Psi$ (thought of as an $\dim V$-dimensional representation of $L_E \times \SL_2(\C)$), then under the global theta lift with respect to $(\chi_V, \chi_W, \psi)$, the global theta lift of $\pi$ (which is a summand in $\mathcal{A}_2(\U(W))$) has A-parameter
 \[  \chi_V \cdot (\chi_W^{-1} \Psi  \oplus  S_{\dim W - \dim V}). \]
  \end{con}
 
 Given Arthur's conjecture, this is largely a local unramified issue.  What is subtle about Adam's conjecture is its local analog (which we don't discuss here)
  \vskip 5pt
  
  \subsection{\bf Examples}
  We shall illustrate Arthur's conjecture with several families of examples in low rank in the rest of the lecture. For the groups $\SO_5$, $\U_3$ and $G_2$, we shall write down a family of nontempered A-paremeters. For each such A-parameter $\psi$, we will examine the consequences of Arthur's conjecture. This will involve determining:
  \vskip 5pt
  
  \begin{itemize}
  \item the local and global component groups associated to $\psi$;
  \item the quadratic character $\epsilon_{\psi}$;
  \item the size of the local A-packets and the structure of the submodule $\mathcal{A}_{\psi}$. 
  \end{itemize}
  We will then see if the description of these A-parameters provide some clues to how the A-packets and $\mathcal{A}_{\psi}$ may be constructed. 
 \vskip 10pt

\subsection{\bf  Saito-Kurokawa example}
Let $G = \SO_5 = \PGSp_4$, so that its Langlands dual group is $G^{\vee}= \Sp_4(\C)$. 
We have the subgroup  
\[ \SL_{2}(\C) \times \SL_{2}(\C)  \subset \Sp_4(\C) = G^{\vee}. \]
These two commuting $\SL_2$'s play symmetrical roles, as they correspond to a pair of orthogonal long roots in the $C_2$ root system.
 
\vskip 5pt

We will  consider A-parameters of the form:
\[  \psi = \rho \times {\rm Id} :  L_k \times \SL_2(\C) \longrightarrow \SL_2(\C) \times \SL_2(\C) \subset G^{\vee} = \Sp_4(\C). \]
Such an A-parameter $\psi$ is specficied by giving an (admissible)  homomorphism 
\[  \rho: L_k \longrightarrow  \SL_2(\C). \]
Note also that
\[  Z_{G^{\vee}}(\psi) =  Z_{\SL_2}(\rho) \times Z_{\SL_2}  \]
Hence, the A-parameter $\psi$ is discrete if and only if $Z_{\SL_2}(\rho)$ is finite, or equivalently if the 2-dimensional representation of $L_k$ afforded by $\rho$ is irreducible. 
By our (basic hypothesis), to give such a $\rho$ is the same as giving a cuspidal representation $\tau  = \tau_{\rho}$ of $\GL_2$ with trivial central character, i.e. a cuspidal representation of $\PGL_2$.

\vskip 5pt

A discrete A-parameter of $G = \SO_5$ of the above form is called a Saito-Kurokawa A-parameter. We have just seen that such A-parameters are parametrized by cuspidal representations of $\PGL_2$.

\vskip 5pt
Given a parameter $\psi = \psi_{\tau}$,  let us compute the various quantities that appear in Arthur's conjecture.
As we saw above
\[  S_{\psi} = (Z_{\SL_2}(\rho_{\tau}) \times Z_{\SL_2})/ Z_{\Sp_4}  = (\mu_2 \times \mu_2)/ \Delta \mu_2 = \mu_2 \]
Likewise  the local component groups $S_ {\psi_{\tau}}$ are given by
\[ S_{\psi_{\tau,v}} =  \begin{cases}
\mu_2, \text{ if $\rho_{\tau_v}$ is irreducible;} \\
1, \text{  if $\rho_{\tau_v}$ is reducible.} \end{cases} \]
The condition that $\rho_{\tau,v}$ be irreducible  is equivalent to
$\tau_v$ being a discrete series representation of $\PGL_2(F_v)$.

\vskip 5pt
\subsubsection{\bf Local Arthur packets.}
Now Arthur's conjecture  predicts that for each place $v$, 
the local $A$-packet $A_{\psi_{\tau,v}}$ has the form:
$$A_{\psi_{\tau,v}}=\left\{ \begin{array}{ll}
  \{\pi_{\tau_v}^+,\pi_{\tau_v}^-\}, & {\rm if \,\,\tau_v\,\, is\,\, discrete \,\, series,} \\
 \{\pi_{\tau_v}^+\}, & {\rm if \,\,\tau_v\,\, is\,\, not \,\, discrete\,\, series}
  \end{array}\right.
$$ 
Here, $\pi_{\tau_v}^+$ is indexed by the trivial character of $S_{\psi_{\tau, v}}$.
\vskip 5pt

Of course, we know what $\pi^+_v$ has to be for almost all $v$: it is irreducible unramified
with Satake parameter $s_{\psi_{\tau, v}}$. This unramifed representation $\pi_v^+$
is the unramfied constituent of the induced representation
 \[  I_P(\tau_v, 1/2)
 =  {\rm Ind}_P^G  |-|_v^{1/2} \otimes \tau_v. \]
 where $P = MN$ is the Siegel parabolic subgroup of $\SO_5$ with Levi factor $\GL_1 \times \SO_3 = \GL_1 \times \PGL_2$. 
 From this, we see that the representations in the global A-packet are nearly 
equivalent to the constitutents of ${\rm Ind}_{P(\A)}^{G(\A)} |-|^{1/2} \otimes \tau$.  Moreover, their local components  are nontempered for almost all $v$.
 
\vskip 15pt

\subsubsection{\bf Global $A$-packets.}
Let $S_{\tau}$ be the set of places $v$ where $\tau_v$ is discrete series,
so that the global $A$-packet has $2^{\# S_{\tau}}$ elements. 
 
\vskip 5pt 
To compute  the 
multiplicity $m_{\eta}$ of $\pi_{\eta} \in A_{\psi_{\tau}}$  
we need to know the quadratic character $\epsilon_{\psi_{\tau}}$ of 
$S_{\psi_{\tau}}$. By a short computation (which you should do),
$\epsilon_{\psi_{\tau}}$ is the non-trivial character of 
$S_{\psi_{\tau}} \cong \mu_2$
if and only if $\epsilon(1/2, \tau) = -1$. Here $\epsilon(s, \tau)$ is the global $\epsilon$-factor of $\tau$.
\vskip 10pt

Now if $\pi = \otimes_v \pi_{\tau_v}^{\epsilon_v} \in  A_{\psi_{\tau}}$, then the 
multiplicity associated to $\pi$ by Arthur's conjecture is: 
\[ m(\pi) = \begin{cases}
1, \text{  if $\epsilon_{\pi}: = \prod_v \epsilon_v = \epsilon(1/2, \tau)$;} \\
0, \text{  if $\epsilon_{\pi} = -\epsilon(1/2, \tau)$.} \end{cases} \]
Thus, we should have:
\[ \mathcal{A}_{\psi_{\tau}} \cong \bigoplus_{\pi \in A_{\psi_{\tau}}:\epsilon_{\pi} = \epsilon(1/2, \tau)} \pi. \]
\vskip 10pt

\subsubsection{\bf Construction}
How can one construct the A-packets $A_{\psi_{\tau}}$ and the space $\mathcal{A}_{\psi_{\tau}}$? It seems that 
we need a lifting to go from $\tau$ to these square-integrable representations of $\PGSp_4$. Since the theta correspondence is 1-to-1, one cannot hope to use theta correspondence to go from $\tau$ to $\mathcal{A}_{\psi_{\tau}}$, never mind the fact that $\PGL_2 \times \SO_5$ is not a dual pair in a symplectic group. 
\vskip 5pt

We need an intermediate step: the Shimura correspondence, or rather its automorphic description by Waldspurger \cite{W1, W3}. 
Using the theta correspondence for $\Mp_2 \times \SO_3$, Waldspurger was able to provide a classification of the constituents of $\mathcal{A}_2(\Mp_2)$  in the style of Arthur's conjecture. 
\vskip 5pt

More precisely,
  $\tau$ gives rise to a packet of cuspidal representations on $\Mp_2$, whose structure is exactly the same as that of the Saito-Kurokawa packets. Namely,  for each place $v$, one has a local packet of irreducible representations of $\Mp_2(k_v)$:
\[  \tilde{A}_{\tau_v} = \begin{cases}
\{ \sigma_{\tau_v}^+, \sigma_{\tau_v}^- \}, \text{  if $\tau_v$ is discrete series;} \\
\{ \sigma_{\tau_v}^+ \}, \text{  if $\tau_v$ is not discrete series.} \end{cases} \]
We call these the Waldspurger packets.
One can form the global packet as a restricted tensor product of the local ones, and one gets a submodule
\[  \tilde{\mathcal{A}}_{\tau} = \bigoplus_{\pi \in \tilde{A}_{\tau}:\epsilon_{\sigma} = \epsilon(1/2, \tau)} \sigma \subset \mathcal{A}_{cusp}(\Mp_2). \]
 . 
\vskip 5pt

Observe the formal similarity between the structure of the Waldspurger packets and the Saito-Kurokawa ones. Given this, and taking note that one has a dual pair
\[  \Mp_2 \times \SO_5 \]
(which is the next step of the $\SO_{2n+1}$ Rallis tower),
it is then not surprising that the local Saito-Kurokawa packets can be constructed as local theta lifts of the local Waldspurger packet: one sets
\[  \pi_{\tau_v}^{\epsilon_v}  = \theta_{\psi_v}(\sigma_{\tau_v}^{\epsilon_v}). \]
These local theta lifts are nonzero because we are in the stable range. Likewise, the Saito-Kurokawa submodule $\mathcal{A}_{\psi_{\tau}}$ can be constructed as the global theta lift of the submodule $\tilde{\mathcal{A}}_{\tau} \subset \mathcal{A}_{cusp}(\Mp_2)$. This was first studied by Piatetski-Shapiro \cite{PS}, but see \cite{G2} for a more refined discussion.
\vskip 10pt 

\subsection{\bf $\U_3$: Howe-PS example}
Now we carry out the same analysis for a family of nontempered A-parameters of $G = \U_3$ (relative to $E/k$) which will explain the Howe-PS example we discussed. 
\vskip 5pt

The Langlands dual group  of $G= \U_3$ is $\GL_3(\C)$, but we need to work with the L-group ${^L}G = \GL_3(\C) \rtimes {\rm Gal}(E/k)$. The A-parameters of $G$ are then 
\[  \psi: L_k \times \SL_2(\C) \longrightarrow  {^L}G. \]
Thankfully, by \cite[Thm. 8.1]{GGP}, the equivalence class of $\psi$ is determined by the equivalence class of its restriction to $L_E$, so we can simply consider
\[  \psi:  L_E \longrightarrow G^{\vee} =\GL_3(\C). \]
In other words, an A-parameter of $G = \U_3$ is simply a 3-dimensional semisimple representation of $L_E$. But this representation needs to satisfy an extra condition: it should be conjugate orthogonal.  In addition, for it to be discrete, $\psi$ should be multiplicity-free.
\vskip 5pt

Clearly, one has a subgroup
\[  \GL_1(\C) \times \GL_2(\C) \subset \GL_3(\C). \]
We are going to build a discrete  A-parameter 
\[  \psi: L_E \times \SL_2(\C) \longrightarrow GL_1(\C) \times \GL_2(\C) \subset \GL_3(\C), \]
so that $\psi(\SL_2(\C)) = \SL_2(\C) \subset \GL_2(\C)$. As a 3-dimensional representation, $\psi$ takes the form
\[  \psi= \psi_{\chi,\mu}  := \mu \oplus \chi \otimes S_2 \]
where $S_2$ denotes the 2-dimensional irreducible representation of $\SL_2(\C)$. The conjugate-orthogonal condition amounts to requiring that
\vskip 5pt

\begin{itemize}
\item $\mu$ is a conjugate-orthogonal 1-dimensional character of $L_E$, which by our (basic hypothesis) corresponds to an automorphic character of $\A_E^{\times}$ trivial on $\A_k^{\times}$;
\item $\chi$ is a conjugate-symplectic 1-dimensional character of $L_E$, which corresponds by our (basic hypothesis)  to an automorphic character of $\A_E^{\times}$ whose restriction to $\A_k^{\times}$ is $\omega_{E/k}$.
\end{itemize}
Thus, such a $\psi = \psi_{\chi, \mu}$ is specified by the pair $(\chi, \mu)$ satisfying the above properties.

\vskip 5pt

\subsubsection{\bf Component groups and A-packets}
The global component group of $\psi$ is
\[  S_{\psi} = \mu_2\]
and the local component groups are:
\[  S_{\psi_v} = \begin{cases}
\mu_2, \text{  if $v$ is inert in $E$;}\\
1, \text{  if $v$ splits in $E$.} \end{cases} \]
So the local A-packets have the form
\[ A_{\psi_v} =  \begin{cases}
\{  \pi_v^+, \pi_v^- \}, \text{  if $v$ is inert in $E$;} \\
\{ \pi_v^+ \}, \text{  if $v$ splits in $E$.} \end{cases}\]
Moreover, for almost all inert places $v$, the representation $\pi_v^+$ is the unramified representation contained in the principal series representation 
\[  {\rm Ind}_B^G \chi |-|^{1/2} \otimes \tilde{\mu}, \]
where $\tilde{\mu}$ is $\mu$ regarded as a character of $E^1$, via the standard isomorphism $E_v^{\times}/ k_v^{\times} \cong E_v^1$.  Observe that the global A-packet $A_{\psi}$  has infinitely many elements in this case.
\vskip 10pt

\subsubsection{\bf Multiplicity formula}
To work out the multiplicity formula, we need to work out the quadratic character $\epsilon_{\psi}$. A short and instructive computation shows that
$\epsilon_{\psi}$ is the trivial character of $S_{\psi} = \mu_2$ if and only if
\[  \epsilon_E(1/2, \chi \mu^{-1}) = 1. \]
So Arthur's conjecture predicts that
\[   \mathcal{A}_{\psi_{\chi,\mu}} \cong \bigoplus_{\pi \in A_{\psi}: \epsilon(\pi) = \epsilon_E(1/2, \chi\mu^{-1})}  \pi. \] 

\vskip 10pt

\subsubsection{\bf Construction via theta lifts}
 Adam's conjecture suggests that
  the theta correspondence for $\U_1 \times \U_3$ that we discussed in Lectures 2 and 3 can be used to construct the local A-packets and the submodule
$\mathcal{A}_{\psi_{\chi,\mu}}$.   
\vskip 5pt

Let us fix our skew-Hermitian space $W$ over $k$ and consider an inert place $v$ of $k$.  Recall that there are two rank 1 Hermitian spaces $V_v^+$ and $V_v^-$ over $k_v$ for such an inert place $v$. Roughly speaking, the local A-packet $A_{\psi_{\chi_v, \mu_v}}$ should be constructed as the local theta lift of a particular character of $\U(V_v^+)$ and $\U(V_v^-)$, under the two theta correspondences for $\U(V_v^+) \times \U(W)$ and $\U(V_v^-) \times \U(W)$.  To make this precise, we need to answer a few questions:

\vskip 5pt

\begin{itemize}
\item Which pair of splitting characters $(\chi_V, \chi_W)$ should we use for the theta correspondence?
\vskip 5pt

\item Having fixed $(\chi_V, \chi_W)$, which character of $\U(V_v^+) = \U(V_v^-)$ should we start with?
\vskip 5pt

\item How should we label the two representations in the local A-packet, i.e. which of these two theta lifts is $\pi_v^+$, so as to achieve the predicted multiplicity formula?
\end{itemize}

Based on our understanding of the theta correspondence from Lectures 2 and 3, can you answer these questions?  
\vskip 10pt

\subsection{\bf Example of $G_2$}
We conclude this section by describing 2 families of $A$-parameters of the split exceptional group $G_2$.

\vskip 5pt

\subsubsection{\bf Some structural facts}
The Langlands dual group of $G$
is $G_2(\mathbb{C})$.  We list a few relevant facts about the structure of $G_2(\C)$, referring to its root system here for justification:

 
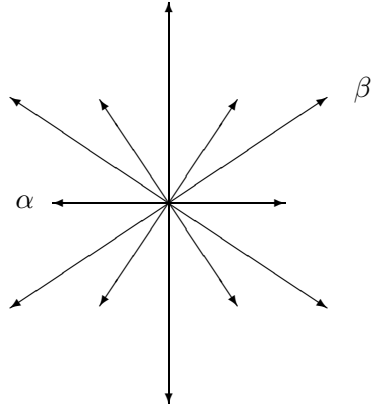
\begin{picture}(300,200)(-100,0)

\put(100,100){\vector(1,0){44}}
\put(100,100){\vector(-1,0){44}}
\put(100,100){\vector(0,1){76}}
\put(100,100){\vector(0,-1){76}}
\put(100,100){\vector(3,2){60}}
\put(100,100){\vector(-3,2){60}}
\put(100,100){\vector(-3,-2){60}}
\put(100,100){\vector(3,-2){60}}
\put(100,100){\vector(2,-3){26}}
\put(100,100){\vector(2,3){26}}
\put(100,100){\vector(-2,3){26}}
\put(100,100){\vector(-2,-3){26}}
\put(42,98){$\alpha$} 
\put(170,140){$\beta$}
\end{picture}

\vskip 5pt

\begin{itemize}
\item The root system of $G_2$ contains a mutually orthogonal pair of long and short roots, giving rise to a commuting pair of $\SL_2$'s (as in the case of $\Sp_4(\C)$)
\[ (\SL_{2,l} \times \SL_{2,s})/ \mu_2^{\Delta}  \subset G_2. \]
The difference is that in the $\Sp_4$ case, the two roots involved have the same length (they are both long), whereas here they are of different length. Hence, these two $\SL_2$'s are not conjugate to each other. Further,  the centralizer of one of these $\SL_2$'s is the other $\SL_2$.
\vskip 5pt

\item The 6 long roots of the $G_2$ root system gives an $A_2$ root system, reflecting the fact that $G_2(\C)$ contains a subgroup $\SL_3(\C)$.
Observe that $\SL_{2,l}(\C)  \times_{\mu_2} T \subset \SL_3(\C)$ (where $T$ is the diagonal torus of $\SL_{2,s}$)  but $\SL_{2,s}$ is not contained in $\SL_3$ (even after conjugation). Moreover, the normalizer of $\SL_3(\C)$ in $G_2(\C)$ contains $\SL_3(\C)$ with index 2. Indeed, an element in the non-identity component is given by the longest Weyl group element of $G_2$. This is also the element $(w,w) \in \SL_{2,l} \times_{\mu_2} \SL_{2,s}$, where $w$ is the standard Weyl element in $\SL_2$.  The conjugation action of this element on $\SL_3(\C)$ is an outer automorphism. Hence one has containment
\[   \SL_{2,l} \times_{\mu_2} N_{\SL_2}(T) \subset N_{G_2}(\SL_3) = \SL_3 \rtimes \Z/2\Z \subset G_2(\C). \]
\vskip 5pt

\item The smallest faithful (irreducible) algebraic representation of  $G_2$ is 7-dimensional; in fact one has
\[  G_2 \hookrightarrow \SO_7. \]
The weights of this 7-dimensional representation are the short roots and the zero vector. 
\vskip 5pt

When restricted to the subgroup $\SL_3$, this 7-dimensional representations decomposes as :
\[  (std_3) \oplus 1 \oplus (std_3)^{\vee} \]
where $(std_3)$ is a 3-dimensional irreducible representation of $\SL_3$ and $(std_3)^{\vee}$ is its dual.  
\vskip 5pt
When restricted to the subgroup $\SL_{2,l} \times_{\mu_2} \SL_{2,s}$, it decomposes as:
\[    (std_2) \otimes (std_2)  \oplus  1 \otimes {\rm Sym}^2(std_2) \]
where $(std_2)$ denotes the 2-dimensional representation of $\SL_2$.
\vskip 5pt
\item consider the adjoint action of $G_2$ on its Lie algebra $\mathfrak{g}_2$. When restricted to the subgroup $\SL_3$,
\[  \mathfrak{g}_2 = \mathfrak{sl}_3 \oplus (std_3) \oplus (std_3)^{\vee}.  \]
When restricted to the subgroup $\SL_{2,l} \times_{\mu_2} \SL_{2,s}$, one has
\[  \mathfrak{g}_2 = \mathfrak{sl}_{2,l} \oplus \mathfrak{sl}_{2,s} \oplus (std_2) \otimes {\rm Sym}^3(std_2). \]
\end{itemize}

\vskip 5pt

\subsubsection{\bf  Some A-parameters}
Now suppose that $\tau$ is a cuspidal representation of $PGL_2$, which by our (basic hypothesis) corresponds to
 an irreducible  representation
\[ \rho_{\tau} : L_F \longrightarrow SL_2(\C). \]
Using $\rho_{\tau}$, we can build 2 different nontempered A-parameters of $G_2$, depending on whether $\SL_2(\C)$ is mapped to $\SL_{2,l}$ or $\SL_{2,s}$:
\[  \psi_{\tau, s}:  L_k \times \SL_2(\C) \longrightarrow \SL_{2,l} \times \SL_{2,s}   \subset G_2(\C) \]
or 
\[  \psi_{\tau,l} :  L_k \times \SL_2(\C) \longrightarrow \SL_{2,s} \times \SL_{2,l}   \subset G_2(\C) .\]
We call $\psi_{\tau, s}$ the short root A-parameter and $\psi_{\tau, l}$ the long root A-parameter associated to $\tau$.
\vskip 5pt

\subsubsection{\bf Short root A-parameters}
Now let's work out the consequences of Arthur's conjecture for the short root A-parameter.
We have:
\[  S_{\psi_{\tau, s}} \cong \mu_2  \]
 and for a place $v$ of $k$,
 \[  S_{\psi_{\tau, s,v}} = \begin{cases}
 \mu_2, \text{  if $\tau_v$ is discrete series (i.e. $\rho_{\tau, v}$ is irreducible);} \\
 1 \text{  if $\tau_v$ is not discrete series (i.e.  $\rho_{\tau,v}$ is reducible).}  \end{cases} \]

 \vskip 5pt
\subsubsection{\bf Local short root A-packets.}
Now Arthur's conjecture  predicts that for each place $v$, 
the local $A$-packet $A_{\tau,s, v}$ has the form:
$$A_{\tau,s, v}=\left\{ \begin{array}{ll}
 \{\pi_{\tau_v}^+,\pi_{\tau_v}^-\}, & {\rm if \,\,\tau_v\,\, is\,\, discrete \,\, series,}\\
 \{\pi_{\tau_v}^+\}, & {\rm if \,\,\tau_v\,\, is\,\, not \,\, discrete\,\, series}
   \end{array}\right.
$$ 
Here, $\pi_{\tau_v}^+$ is indexed by the trivial character of $S_{\tau,v}$.
Moreover,  we know what $\pi^+_v$ has to be for almost all $v$. Indeed, for almost all $v$ where $\tau_v$ is unramified, $\pi_v^+$
is the unramified representation with Satake parameter $s_{\psi_{\tau,v}}$, 
and this representation is the unramified constituent of 
\[ 
I_P(\tau_v,1/2) = {\rm Ind}_P^{G_2} \tau_v \cdot |\det|^{1/2} \]
where $P$ is the Heisenberg parabolic subgroup of $G_2$ with Levi factor $\GL_2$. 
 \vskip 5pt

\subsubsection{\bf Global short root $A$-packets.}
Let $S_{\tau}$ be the set of places $v$ where $\tau_v$ is discrete series,
so that the global $A$-packet has $2^{\# S_{\tau}}$ elements. To describe 
the 
multiplicity of $\pi_{\eta} \in A_{\tau, s}$ in $L^2_{\psi_{\tau}}$, 
we need to know the quadratic character $\epsilon_{\psi_{\tau, s}}$ of 
$S_{\psi_{\tau, s}}$. It turns out that
$\epsilon_{\psi_{\tau, s}}$ is the non-trivial character of 
$S_{\psi_{\tau}} \cong \mu_2$
if and only if $\epsilon(1/2, \tau) = -1$. 
\vskip 5pt

Now if $\pi = \otimes_v \pi_v^{\epsilon_v} \in  A_{\tau,s}$, then the 
multiplicity associated to $\pi$ by Arthur's conjecture is: 
\[ m(\pi) = \begin{cases}
1, \text{  if $\epsilon_{\pi}: = \prod_v \epsilon_v = \epsilon(1/2, \tau)$;} \\
0, \text{  if $\epsilon_{\pi} = -\epsilon(1/2, \tau)$.} \end{cases} \]
Thus, Arthur's conjecture predicts that:
\[\mathcal{A}_{\psi_{\tau,s}} \cong  \bigoplus_{\pi \in A_{\tau}:\epsilon_{\pi} = \epsilon(1/2, \tau)} \pi. \]

\vskip 5pt

\subsubsection{\bf Construction of short root A-packets}
Observe that the structure of these $A$-packets  is thus entirely the same as that of the Saito-Kurakawa packets for $\SO_5$. Since the Saito-Kurokawa packets were constructed as theta liftings of Waldspurger's packets on $\Mp_2$, one might guess that one can construct the short root A-packets of $G_2$ by lifting from the corresponding packets on $\Mp_2$. 
But is $\Mp_2 \times G_2$ a reductive dual pair?
\vskip 5pt

Well, it turns out that one may consider the dual pair 
\[  \Mp_2 \times \O_7.  \]
Recalling that $G_2 \hookrightarrow \SO_7$, we may consider theta lifts from $\Mp_2$ to $\O_7$, followed by restriction of representations from $\O_7$ to $G_2$. Somewhat amazingly, this restriction does not lose much information. In other words, one may consider the commuting pair
\[  \Mp_2 \times  G_2 \]
and restrict the Weil representation of $\Mp_2 \times \O_7$ to it. Such a construction was first conceived by Rallis and Schiffman, but the full analysis was completed in \cite{GG}. In this way, it was shown in \cite{GG} that one may construct the A-packets and the corresponding submodule in $\mathcal{A}_2(G_2)$.
\vskip 10pt

\subsubsection{\bf Long root A-parameters}
The main project for this course is the analysis and construction of the long root A-packets of $G_2$.
In particular, the first task of the project is to work out the prediction of Arthur's conjecture  for the long root A-parameter $\psi_{\tau,l}$, and then specialize to the case when $\tau$ is dihedral.
\vskip 5pt

We list the expected answers here, leaving it as a series of guided exercises:
\vskip 5pt

\begin{itemize}
\item the global and local component groups are the same as for the short root A-parameters; so the local A-packets have 2 or 1 elements depending on whether $\tau_v$ is discrete series or not.

\item the quadratic character $\epsilon_{\psi_{\tau,l}}$ is trivial if and only if 
\[ \epsilon(1/2, \tau, {\rm Sym}^3) = \epsilon(1/2, {\rm Sym}^3 \rho_{\tau}) =1. \]
So we see that the ${\rm Sym}^3$-epsilon factor appears in the Arthur multiplicity formula.
\end{itemize}

\vskip 5pt

\subsection{\bf Dihedral long root A-parameters}
We  now suppose that $\tau$ is a dihedral cuspidal representation relative to a quadratic field extension $E/k$. This can be interpreted in one of the following equvialent ways:
\begin{itemize}
\item $\rho_{\tau} \cong {\rm Ind}_{W_E}^{W_k} \chi$ for some 1-dimensional character $\chi$ of the global Weil group $W_E$ (which is supposedly a quotient of $L_E$).

\item $\tau \otimes \omega_{E/k} \cong \tau$.
\end{itemize}
The fact that $\det \rho_{\tau} = 1$ implies that  when regarded as an automorphic character of $\A_E^{\times}$, $\chi|_{\A_k^{\times}} = \omega_{E/k}$, i.e. $\chi$ is a conjugate-symplectic automorphic character. The image of $\rho_{\tau}$ is contained in the normalizer $N_{\SL_2}(T)$, where $T$ is a maximal torus of $\SL_2$. 
Now we observe:
\vskip 5pt

\begin{itemize}
\item When $\tau$ is dihedral as above, the long root A-parameter $\psi_{\tau, l}$ factors as:
\[  \psi_{\tau,l} : L_k \twoheadrightarrow W_k \longrightarrow N_{\SL_{2,s}}(T) \times_{\mu_2} \SL_{2,l}(\C)  \subset \SL_3(\C) \rtimes \Z/2\Z \subset G_2(\C). \]
This follows from one of the structural facts we recall about $G_2(\C)$.
\vskip 5pt

\item $\SL_3(\C) \rtimes \Z/2\Z \subset  \GL_3(\C) \rtimes \Z/2\Z = {^L} \U_3$.

\item Hence the long root A-parameter $\psi_{\tau,l}$ of $G_2$ factors through the L-group of $\U_3$, thus giving rise to an A-parameter for $\U_3$. Moreover, when restricted to $W_E$, one obtains a 3-dimensional representation of $W_E \times \SL_2(\C)$ of the form
\[    \psi_{\tau, l}|_{W_E} = \chi^{-2} \oplus  \chi \otimes S_2. \]
In other words, one obtains  a Howe-PS A-parameter for $\U_3$. 
\end{itemize}

\vskip 10pt

Said in another way, one could start with a Howe-PS A-parameter $\psi$ for $\U_3$ with $\psi(L_E) \subset \SL_3(\C)$ (or equivalently, giving rise to representations of $\PU_3$). 
The composition of $\psi$ with the natural inclusion
\[  \SL_3(\C) \rtimes \Z/2\Z = N_{G_2}(\SL_3) \hookrightarrow G_2(\C) \]
then gives a long root A-parameter whose associated $\tau$ is dihedral with respect to $E/k$. 

\vskip 5pt

This suggests the following question: 
\vskip 5pt

\noindent{\bf Question:} Is it possible to construct the local and global long root A-packets of $G_2$ by lifting from the Howe-PS packets for $\U_3$, and then verify the Arthur multiplicity formula?
\vskip 5pt
Addressing this question is the main project for this course.

 \vskip 15pt
 
 \section{\bf More Exercises}
In this final section, we give a list of extended exercises.
\vskip 10pt

 \subsection{\bf Lecture 1}  $\text{}$
\vskip 5pt

\noindent (1) This exercise concerns the parametrization of unramified representations of $G(F)$ in terms of semisimple classes in the Langlands dual group $G^{\vee}$.
\vskip 5pt

(i) Let $T$ be a split torus over a p-adic field $F$ and let $T(\mathcal{O}_F) \subset T(F)$ be the maximal compact subgroup of $T(F)$. 
A character $\chi: T(F) \rightarrow \C^{\times}$ is unramified if $\chi$ is trivial on $T(\mathcal{O}_F)$.
The set of unramified characters of $T(F)$  is thus 
\[  \Hom(T(F)/T(\mathcal{O}_F), \C^{\times}). \]
The dual torus of $T$ is the complex torus defined by
\[  T^{\vee} :=  X^*(T) \otimes_{\Z}  \C^{\times},\]
where $X^*(T) = \Hom(T, \mathbb{G}_m)$. Construct a natural bijection
\[  \Hom(T(F)/T(\mathcal{O}_F), \C^{\times}) \cong T^{\vee}. \]
\vskip 5pt

(ii) Based on (i) and Proposition 1.1 of the lecture notes, deduce that unramified representations of $G= \GL_n(F)$ are paramertrized by semisimple conjugacy classes in $G^{\vee}= \GL_n(\C)$.

\vskip 15pt
\noindent (2) This exercise gives you a chance to work with  unitary groups in low rank.
\vskip 5pt

(i)  In  \S 1.9 of the lecture notes, we gave as examples the quasi-split unitary groups $\U(V^+)$ with $\dim V^+ = 2$ and $3$ and wrote down certain elements as matrices.  Let $B \subset \U(V^+)$ be the  upper triangular Borel subgroup. Compute the modulus character $\delta_B$ as a character of the diagonal torus $T$.

\vskip 5pt

(ii) At the end of \S 1.9 of the lecture notes, we introduced the element $u(x,z)$ as a matrix, but one particular entry of the matrix is given as $\ast$,  as it is determined by the others. Determine the entry $\ast$ explicitly.
\vskip 5pt

(iii) In \S 1.9 of the lecture notes, we described some elements of $\U(V)$ where $V$ is a 3-dimensional Hermitian space.
In fact, from Lecture 2 onwards, we will be working with 3-dimensional skew-Hermitian spaces $W$. As what we did for the Hermitian case, write down elements in the Borel subgroup $B = TU$  of $\U(W)$, with respect to a Witt basis of $W$, i.e. a basis $\{e, w_0, e^*\}$ with $\langle e,e^* \rangle = 1$,  $\langle w_0, w_0 \rangle = \delta$ (with $\delta$ trace $0$) and $\langle e,w_0 \rangle= 0 = \langle e^*, w_0 \rangle$. 
\vskip 5pt

 \vskip 5pt
 \subsection{\bf Lecture 2} $\text{}$
 \vskip 5pt
 
\noindent  (3)  This exercise is based on \S 2.3 which introduces the Heisenberg group $H(W)$ associated to a symplectic vector space $W$.
 \vskip 5pt

(i) Let $W$ be a 3-dimensional skew-Hermitian space as in Problem 2(iii) above, with isometry group $\U(W)$ containing the Borel subgroup $B = TU$. Write down an isomorphism of $U$ with the Heisenberg group associated to a 2-dimensional symplectic space. 
\vskip 5pt

(ii) In \S 2.3, we introduced the representation 
\[  \omega_{\psi} = {\rm ind}_{H(X)}^{H(W)} \psi \]
of a Heisenberg group $H(W)$ on the space $\mathcal{S}(Y)$ of Schwartz functions on $Y$. Prove that this representation is indeed irreducible. 
\vskip 5pt

(iii) In the context of \S 2.3, suppose that $W = W_1\oplus W_2$ is the sum of two smaller symplectic spaces, construct a natural surjective group homomorphism 
\[ f:  H(W_1) \times H(W_2) \longrightarrow H(W). \]
What is the kernel of your homomorphism $f$?
\vskip 5pt

(iv) Let $\omega_{W, \psi}$ be the irreducible representation of $H(W)$ with central character $\psi$. Show that the pullback $f^*(\omega_{W,\psi})$ is isomorphic to $\omega_{W_1,\psi} \otimes \omega_{W_2,\psi}$. 
\vskip 5pt

(v)  One has a natural embedding
\[  f: \Sp(W_1) \times \Sp(W_2) \hookrightarrow \Sp(W). \]
 Deduce that $f^*$ induces a natural isomoprhism of projective representations: 
 \[  A_{W,\psi} \circ f \cong  A_{W_1, \psi} \otimes A_{W_2,\psi} \]
 where $A_{W, \psi}: \Sp(W) \rightarrow \GL(\mathcal{S})/ S^1$ is as constructed in \S 2.3. 
 \vskip 10pt
 
 Indeed, $f$ can be lifted to 
 \[  \tilde{f}: \Mp(W_1) \times \Mp(W_2) \longrightarrow \Mp(W) \]
 so that one has an isomoprhism 
 \[  \tilde{f}^*(\omega_{W,\psi}) \cong \omega_{W_1, \psi} \otimes \omega_{W_2, \psi_2} \]
 of representations of $\Mp(W_1) \times \Mp(W_2)$.

\vskip 15pt

\noindent (4) This exercise asks you to reconcile the different ways that $\Theta(\pi)$ has been presented in the lecture notes.
 \vskip 5pt
 
 In \S 2.5, we have given two descriptions of $\Theta(\pi)$:
 \vskip 5pt
 
 \begin{itemize}
 \item 
 \[ \pi \boxtimes  \Theta(\pi)  = \Omega/ \bigcap_{f \in \Hom_{\U(V)}(\Omega, \pi)} {\rm Ker}(f). \]

\item  $\Theta(\pi) = (\Omega \otimes \pi^{\vee})_{\U(V)}$. 
\end{itemize}
Show that these are equivalent, and prove the ``universal property":
\[  \Hom_{\U(V) \times \U(W)}(\Omega, \pi \otimes \sigma) \cong \Hom_{\U(W)}(\Theta(\pi),\sigma) \]
for any smooth representation $\sigma$ of $\U(W)$. 

\vskip 15pt

\noindent (5) Deduce Corollary 2.5 (dichotomy) from Theorem 2.4 (conservation relation).
\vskip 15pt

\noindent (6) This problem introduces the Schrodinger model of the Weil representation of a dual pair $\U(V) \times \U(W)$, when one of the spaces is split of even dimension.
\vskip 15pt

In \S 2.4, we wrote down some formulas in the Schrodinger model for the Weil representation $\omega_{\psi}$ of a metaplectic group; this model is based on a Witt decomposition of the symplectic space. In \S 2.5, we considered a dual pair $\U(V) \times \U(W)$ with its splitting
\[  \tilde{\iota}: \U(V) \times \U(W) \longrightarrow \Mp(V \otimes_E W)\]
associated to a pair $(\chi_V, \chi_W)$.
One can ask if we can inherit the formulas in the Schrodinger model and write down the action of some elements of $\U(V) \times \U(W)$. For this to be possible, there must be some compatibility between the Witt decomposition we used on $V \otimes W$ and the map $\iota: \U(V) \times \U(W) \rightarrow \Sp(V \otimes_E W)$. 
\vskip 5pt

More precisely, suppose $V$ is a split Hermitian space and we fix a Witt decomposition $V = X \oplus Y$. Then we inherit a Witt decomposition of $V \otimes W$:
\[  V \otimes_E W = (X \otimes_E W) \oplus (Y \otimes_E W). \]
Relative to this Witt decomposition, the Schrodinger model of the Weil representation is realized on $\mathcal{S}(Y \otimes W)$ and one can write down explicit formulas for the elements of 
\[  \U(W) \times P(X) \subset P(X \otimes W), \]
where $P(X)$ is the Siegel parabolic subgroup of $\U(V)$ stabilzing $X$.
\vskip 5pt

Consider the case when $W =  E \cdot w= \langle \delta  \rangle$  (with $Tr(\delta) =0$) is 1-dimensional and $V = E e \oplus E e^*$ is the split skew-Hermitian space of dimension $2$.  From the formulas of the Schrodinger model, deduce (as much as you can) the following actions of  $\U(W) \times B(E\cdot e))$ on $\mathcal{S}(Y \otimes V) = S(E e^* \otimes w)$  (relative to the fixed $(\chi_V, \chi_W)$):
\begin{itemize}
\item for $g \in \U(W) = E^1$,
\[  (g \cdot f)( x)  = \chi_V( i(g)) f( g^{-1} x), \]
where $i: E^1 \cong E^{\times}/F^{\times}$ is the inverse of the isomorphism $i^{-1}: x \mapsto x/x^c$. 
\vskip 5pt

\item For $t(a) \in T$, $a \in E^{\times}$,
\[  (t(a) \cdot f) (x) =  \chi_W(a) \cdot  |a|_E^{1/2} \cdot f(a^c \cdot x). \]
\vskip 5pt

\item For $u(z) \in U$, with $z \in E$ and  $Tr(z) = 0$,
\[  (u(z) \cdot f)(x) = \psi( \delta \cdot z\cdot N(x)) \cdot f(x). \]
\end{itemize}
To be honest, since we did not explicate the definition of the splitting associated to $(\chi_V, \chi_W)$, you could not really show the above formulas in full , but you can at least deduce those parts of the formula without $\chi_V$ or $\chi_W$.  The effects of the choice of $(\chi_V, \chi_W)$ can be seen from the first two formulas. 
\vskip 5pt

Actually, this exercise is setting the scene for Problem 7 below, so you may take the formulas above as  given.
\vskip 15pt

\noindent (7) This exercise continues from Problem 6. It is the first exercise that allows you to work with the Weil representation and to calculate some theta lifts.
\vskip 5pt

Consider the dual pair $\U(V) \times \U(W)$ as in Problem 6, so that $W= E \cdot w =\langle \delta \rangle$ and $V$ the split 2-dimensional Hermitian space. Fix a pair of splitting characters $(\chi_V, \chi_W)$ and consider the associated Weil representation $\Omega = \Omega_{\chi_V, \chi_W, \psi}$ of $\U(V) \times \U(W)$. Because $\U(W)$ is compact, $\Omega$ is semisimple as a $\U(W)$-module and we can write:
\[  \Omega = \bigoplus_{\mu \in {\rm Irr}(\U(W))} \mu \otimes \Theta(\mu). \]
The goal is to understand  $\Theta(\mu)$  as much as possible.
\vskip 5pt

Using the formulas for the Weil representation $\Omega = \Omega_{\chi_V, \chi_W, \psi}$ from Problem 6 above, 
\vskip 5pt

(i) Compute $\Omega_U$ (the $U$-coinvariants of $\Omega$) as a module for $\U(W) \times T$. 

\vskip 5pt

(ii) For any nontrivial character $\psi'$ of $U \cong F \cdot \delta^{-1}$, compute $\Omega_{U, \psi'}$ as a module for $\U(W) \times Z(\U(V))$. (Note that there are two orbits of such nontrivial characters $\psi'$ under the conjugation action of $T(F)$). 
\vskip 5pt

(iii) Using your answers from (i) and (ii), show that $\Theta(\mu)$ is nonzero irreducible for any irreducible character $\mu$ of $\U(W) = E^1$. Moreover, show that $\Theta(\mu)$ is supercuspidal if and only if $\mu\ne \chi_V \circ i$ (see Problem 6 for the definition of $i$).  

\vskip 5pt

(iv) 
Show that
\[  \Theta(\chi_V \circ  i) \hookrightarrow I(\chi_W) := {\rm Ind}_{B(E \cdot e)}^{\U(V)} \chi_W. \]
Indeed, your proof should suggest an explicit description of this embedding. More precisely, show that the natural map
\[  f  \mapsto  \left( h \mapsto (\Omega(h) f)(0) \right) \]
gives a nonzero equivariant map
\[  \Omega \longrightarrow (\chi_V \circ i) \otimes I(\chi_W),  \]
 thus inducing the embedding of $\Theta(\chi_V \circ i)$ into $I(\chi_W)$. 
 
\vskip 5pt
The results of this exercise will be used in the next exercise.

\vskip 15pt

\noindent (8) The purpose of this exercise is to indicate a proof of  the Howe duality conjecture and Theorem 2.6 for the dual pair $\U_1 \times \U_1$. It is very  long, but is the highlight of the problem sheet!
As mentioned in the notes, the proof makes use of the doubling seesaw argument (among other things). We will introduce some of these notions in turn.
\vskip 5pt

\begin{itemize}
\item (Seesaw pairs)
Suppose a group $E$ contains two dual pairs $G_1 \times H_1$ and $G_2 \times H_2$ (so $G_i$ is the centralizer of $H_i$ in $E$ and vice versa). Suppose that $G_1 \subset G_2$. Then it follows that $H_1 \supset H_2$. In this situation, we say that the two dual pairs form a seesaw pair, and we often represent this in the following seesaw diagram:

\vskip 5pt

 \[
 \xymatrix{
  G_2  \ar@{-}[dr] \ar@{-}[d] & H_1  \ar@{-}[d] \\
  G_1  \ar@{-}[ur] & H_2}
\]  

\vskip 5pt

In this diagram, the diagonal line represents a dual pair, and the vertical line denotes containment, with the group at the bottom contained in the group at the top. 

\vskip 5pt

\item (Standard example) Here is the standard example of constructing seesaw pairs in the symplectic groups. Suppose that $V_1 \oplus V_2$ is the orthogonal sum of two Hermtian spaces. 
Set $\mathcal{W} = (V_1 \oplus V_2) \otimes W$ (a symplectic space over $F$), and note that
\[ \mathcal{W} = (V_1 \otimes W)  \oplus  (V_2 \otimes W). \]
This gives the following two dual pairs in $\Sp(\mathcal{W})$:
\[  \U(V_1 + V_2) \times \U(W^{\Delta}) \quad \text{and} \quad (\U(V_1) \times \U(V_2))) \times (\U(W) \times \U(W)).  \]
Convince yourself that these form a seesaw pair and draw the relevant seesaw diagram.

\vskip 5pt
\item (Seesaw identity) Suppose one has a seesaw diagram as in the abstract situation above, and let $\Omega$ be a representation of $E$, which we may restrict to $G_1 \times H_1$ and $G_2 \times H_2$. Deduce the following seesaw identity: for $\pi \in {\rm Irr}(G_1)$ and $\sigma \in {\rm Irr}(H_2)$, one has natural isomorphisms
\[  \Hom_{G_1}(\Theta(\sigma), \pi) \cong  \Hom_{G_1 \times H_2}(\Omega, \pi \otimes \sigma)  \cong \Hom_{H_2}(\Theta(\pi), \sigma). \]
Note that in the above identity, $\Theta(\sigma)$ is a representation of $G_2$, whereas $\Theta(\pi)$ is a representation of $H_1$.
This seesaw identity allows one to transfer a restriction problem from one side of the seesaw to the other.

\vskip 5pt

\item (Compatible splittings)  To apply this seesaw identity to the standard example, there is an extra step, because we need to consider splittings of the dual pair into the metaplectic group $\Mp(\mathcal{W})$. To have a splitting of the metaplectic cover for the dual pair $\U(V_1) \times \U(W)$, we need to fix a pair $(\chi_{V_1}, \chi_W)$; likewise, we need to fix $(\chi_{V_2}, \chi'_W)$ for $\U(V_2) \times \U(W)$. Similarly, for the dual pair $\U(V_1 + V_2) \times \U(W^{\Delta})$, we may fix $(\chi_{V_1 + V_2}, \chi_{W^{\Delta}})$. 
So we see that we have 6 splitting characters to fix here. If we were to choose these randomly, then there is no reason for the resulting splittings to be compatbile with each other.
\vskip 5pt

What does being compatible with each other mean? From the viewpoint of the seesaw diagram,  in order to have the seesaw identity, we need to ensure that when the splittings of  a group at the top of the diagram is restricted to the subgroup below it, the restriction agrees with the splitting below. This is to ensure that we still have a seesaw situation in $\Mp(\mathcal{W})$. 
\vskip 5pt

  So for example, we fix the character $\chi_{W^{\Delta}}$ which determines the splitting of $\U(V_1+ V_2)$. When restricted to $\U(V_1) \times \U(V_2)$, the resulting splitting over $\U(V_1)$ and $\U(V_2)$ are both associated with $\chi_{W^{\Delta}}$. This forces us to take
\[  \chi_{W^{\Delta}} = \chi_W = \chi'_W, \]
and we shall denote this by $\chi_W$. Likewise, we fix the characters $\chi_{V_1}$ and $\chi_{V_2}$ which determine a splitting over $\U(W) \times \U(W)$. When restricted to the diagonal $U(W^{\Delta})$, the resulting splitting of the latter is associated with the character $\chi_{V_1} \chi_{V_2}$.  This forces us to take
\[  \chi_{V_1+ V_2} = \chi_{V_1} \cdot \chi_{V_2}. \]

\end{itemize}
\vskip 5pt

(i)  (Doubling seesaw)  Now we   apply the above to the following concrete situation. We place ourselves in the setting of Theorem 2.6. Hence, let $W = \langle b \cdot \delta \rangle$ be a 1-dimensional skew-Hermitian space, and $V = \langle a \rangle$  a 1-dimensional Hermitian space, so that $\epsilon(W) = \omega_{E/F}(b)$ and $\epsilon(V) = \omega_{E/F}(a)$. 
Let   $V^- = \langle -a \rangle$ and  apply the seesaw construction above with $W$ as given,
\[  V_1 =   V, \quad \text{and} \quad V_2 = V^-   \]
so that
\[  V^{\square} := V_1 \oplus V_2  = V \oplus V^- \]
is a split 2-dimensional Hermitian space. We thus have the seesaw diagram:
\vskip 5pt

 \[
 \xymatrix{
  \U(V^{\square})  \ar@{-}[dr] \ar@{-}[d] & \U(W) \times \U(W)   \ar@{-}[d] \\
  \U(V) \times \U(V^-)  \ar@{-}[ur] & \U(W)^{\Delta}}
\]

\vskip 5pt

This is called the {\em doubling seesaw}, because we have doubled $V$ (to yield $V^{\square})$, but note that we have introduced a negative sign in the second copy of $V$, so that the doubled-space $V^{\square}$ is split! Indeed, the diagonally embedded $V^{\Delta}$ is a maximal isotropic subspace, and one has a Witt decomposition
\[  V^{\square} = V^{\Delta} \oplus V^{\nabla} \]
where $V^{\nabla} = \{ (v,-v): v \in V \}$.

\vskip 5pt

(ii) (Doubling seesaw identity)  Choose splitting characters $\chi_V$, $\chi_{V^-}$ and $\chi_W$ as explained above.  In fact, we insist further (as we may) that
\[  \chi_{V} = \chi_{V'} = \chi_W = \gamma \quad \text{(a conjugate-symplectic character of $E^{\times}$).} \]
  Fix a $\chi \otimes \chi' \in {\rm Irr}(\U(V) \times \U(V^-))$ and the character $\chi_V|_{E^1}$ of $\U(W^{\Delta})$. Write down what the doubling seesaw identity gives.

\vskip 5pt

(iii)  (Duality) We have the two decompositions
\[  \Omega_{V, W, \gamma,\psi} = \bigoplus_{\chi} \chi \otimes \Theta_{V,W}(\chi),   \]
and
\[  \Omega_{V^-, W,\gamma, \psi} = \bigoplus_{\chi'}\chi'  \otimes  \Theta_{V^-,W}(\chi').      \]
Express the Weil representation $\Omega_{V^-, W, \gamma, \psi}$  in terms of $\Omega_{V,W, \gamma, \psi}$ and deduce that
\[  \Theta_{V^-, W}(\chi^{-1} \chi_W|_{E^1}) \cong \overline{\Theta_{V,W}(\chi)} \cdot \chi_V|_{E^1}. \] 
Using this, show that  
\[ \Omega_{V, W, \gamma,\psi}  \otimes  \Omega_{V^-, W,\gamma, \psi} =  \bigoplus_{\chi, \chi'}  \chi \otimes {\chi'}^{-1} \chi_W|_{E^1}  \otimes \Hom_{\U(W)}(  \Theta_{V,W}(\chi) \otimes \overline{ \Theta_{V,W}(\chi')}, \C) \]
as a module for $\U(V) \times \U(V^-) \times \U(W^{\Delta})$. 
(Note that here and below, we could have replaced $\chi_V$ and $\chi_W$ by $\gamma$, but we have refrained from doing so, in order to make the dependence of $(\chi_V, \chi_W)$ more transparent).
 \vskip 5pt

(iv) (Siegel-Weil)  
The LHS of the doubling seesaw is 
\[  \Hom_{\U(V) \times \U(V^-)} (\Theta_{W^{\Delta}, V^{\square}}(\chi_V|_{E^1}), \chi \otimes {\chi'}^{-1} \chi_W|_{E^1}). \]
To address this problem, we first need to understand $\Theta_{W^{\Delta}, V^{\square}}(\chi_V|_{E^1})$. This is  the local version of the Siegel-Weil formula and is where Problem 7 comes in.
\vskip 5pt
Using your results in Problem 7(iv), show that 
\[  \Theta_{W^{\Delta}, V^{\square}}(\chi_V|_{E^1}) \hookrightarrow {\rm Ind}_{B(V^{\Delta})}^{\U(V^{\square})} \chi_W =: I(\chi_W) \]

(v) (Principal series)
Show that the principal series $I(\chi_W)$  is reducible and in fact is the direct sum of two irreducible summands. How does one distinguish between those two summands? Which of these two summand is $\Theta_{W^{\Delta}, V^{\square}}(\chi_V|_{E^1})$ equal to? Show that $ \Theta_{W^{\Delta}, V^{\square}}(\chi_V|_{E^1})$ is the unique summand of the induced representation whose $(U, \psi')$-coinvariant is nonzero for $\psi'( z \delta^{-1}) = \psi (b z)$. 
\vskip 5pt

This shows that if $W$ and $W'$ are the two 1-dimensional skew-Hermitian spaces, then 
\[  I(\chi_W) = \Theta_{W, V^{\square}}(\chi_V|_{E^1}) \oplus \Theta_{W', V^{\square}}(\chi_V|_{E^1}). \]

\vskip 5pt

(vi) (Mackey theory) The previous part implies that it will be necessary to understand $I(\chi_W)$ as a module for $\U(V) \times \U(V^-)$. 
This  can be approached by Mackey thoery, as it involves the restriction  of an induced representation. 
\vskip 5pt

Show that $\U(V) \times \U(V^-)$ acts transitively on the flag variety $B(V^{\Delta}) \backslash \U(V^{\square})$ with stabilizer of the identity coset given by $\U(V)^{\Delta}$. From this,  
deduce that as a $\U(V) \times \U(V^-)$-module, $I(\chi_W)$ is isomorphic to 
\[  C^{\infty}_c(\U(V)) \otimes (1 \otimes \chi_W|_{E^1}), \]
i.e. a twist of the regular representation.  Hence, for $\chi \otimes {\chi'}^{-1} \in {\rm Irr}(\U(V) \times \U(V^-))$, one has
\[ \Hom_{\U(V) \times \U(V^-)} (I(\chi_W), \chi \otimes {\chi'}^{-1} \cdot \chi_W|_{E^1}) = \begin{cases}
\C, \text{  if $\chi' =\chi$;} \\
0, \text{  otherwise.}  \end{cases} \]
Moreover, as a $\U(V) \times \U(V^-)$-module, $ \Theta_{W^{\Delta}, V^{\square}}(\chi_V|_{E^1})$ is a submodule of the above twisted regular representation.

\vskip 5pt

(vii) (Howe duality) Using (iii) and (vi), show the following:
\begin{itemize}
 \item For each $\chi \in {\rm Irr}(\U(V))$, $\Theta_{V,W}(\chi)$ is irreducible or $0$; in particular it is either $\chi$ or $0$. 

\item If $\chi \ne \chi'$, then $\Theta_{V,W}(\chi)$  and $\Theta_{V,W}(\chi')$ are disjoint. 
\end{itemize}
This is the Howe duality theorem for $\U(V) \times \U(W)$.
\vskip 10pt

(viii) (Doubling zeta integral)
The remaining issue is to decide for which $\chi$ is $\Theta_{V,W}(\chi) \ne 0$. This requires the use of the doubling zeta integral. In this context, the doubling zeta integral is an explicit integral which defines a nonzero element of 
\[  \Hom_{\U(V) \times \U(V^-)}(I(\chi_W), \chi \otimes \chi^{-1} \chi_W|_{E^1}) = \Hom_{\U(V) \times \U(V^-)}(C^{\infty}_c(\U(V)), \chi \otimes \chi^{-1}). \]
More precisely, we define
\[  Z(s, \chi): I(s, \chi_W) = {\rm Ind}_{B(V^{\Delta})}^{\U(V^{\square})} \chi_W \cdot | -|_E^s  \longrightarrow  \C \]
by
\[  Z(s, \chi)(f_s) = \int_{\U(V)}  f_s( h,1)  \cdot  \overline{\chi(h)}  \, dh. \]
Verify that the integral converges absolutely and defines a nonzero functional
\[  Z(s,\chi) \in \Hom_{\U(V) \times \U(V^-)}(I(s, \chi_W), \chi \otimes \chi^{-1} \chi_W|_{E^1}). \]
Deduce also that $\theta_{V,W}(\chi)$ is nonzero if and only if  $Z(0, \chi)$ is nonzero on $\Theta_{W^{\Delta}, V^{\square}}(\chi_V|_{E^1})$. 
\vskip 10pt

(ix) (Functional equation)
There is a standard $\U(V^{\square})$-intertwining operator
\[  M(s):  I(s, \chi_W) \longrightarrow I(-s, \chi_W),\]
whose precise definition need not concern us here.   There is a normalization of this intertwining operator (which we will not go into here) with the following properties:
\vskip 5pt

\begin{itemize}
\item one has
\[  M(-s) \circ M(s)  = 1. \]
\vskip 5pt

\item At  $s = 0$ (where $I(\chi_W)$ is the sum of two irreducible summands), $M(s)$ is holomorphic and $M(0)$ acts as $+1$ on $\Theta_{W^{\Delta}, V^{\square}}(\chi_V|_{E^1})$ and as $-1$ on the other summand. 
\end{itemize}
\vskip 5pt
A basic result in the theory of the doubling zeta integral is that there is a functional equation:
\[  \frac{Z(-s, \chi^{-1}) (M_s(f_s))}{L_E(\frac{1}{2}-s,  \chi_E^{-1} \chi_W) } =  \epsilon_E(\frac{1}{2}+s, \chi_E \chi_W^{-1}, \psi) \cdot \frac{Z(s, \chi)(f_s)}{L_E(\frac{1}{2}+s,\chi_E \chi_W^{-1})}, \]
 where we recall that $\chi_E(x) = \chi(x / x^c)$ for $x \in E^{\times}$.
 We shall take this as a given.
\vskip 5pt

(x) (Proof of Theorem 2.6)
Using the functional equation and the properties of $M(0)$ recalled in (ix), as well as the relevant results in earlier parts, prove Theorem 2.6.

\vskip 15pt

\noindent (9i) Consider the dual pair $\U(V) \times \U(W)$ as in \S 2.11. Over there, a model for the Weil representation is written down. Understand how the formulas there are deduced from 
Problems 3(v) and  4. 
\vskip 5pt

(ii) Do the exercise formulated at the end of \S 2.11.

\vskip 10pt

\subsection{\bf Lecture 3}  $\text{}$
\vskip 5pt

\noindent (10) The purpose of this exercise is to do the global analog of Problem 7.
Hence, we are considering the dual pair $\U(W) \times \U(V)$ over a number field $k$, with $W = E \cdot w = \langle \delta \rangle$ a 1-dimensional skew-Hermitian space (so $\delta \in E^{\times}$ is a trace 0 element) and $V = E e \oplus Ee^* = X \oplus Y$ is a split 2-dimensional Hermitian space.  As in Problem 6, the global Weil representation $\Omega$ is realized on 
$\mathcal{S}( Y_{\A} \otimes W_{\A}) = \mathcal{S}(\A_E e^* \otimes w)$. The automorphic realization 
\[  \theta: \mathcal{S}( Y_{\A} \otimes W_{\A}) \longrightarrow \mathcal{A}(\U(V) \times \U(W)) \]
of $\Omega$ is defined by
\[  \theta(\phi)(g,h) = \sum_{v \in V_k} (\Omega(g,h) \phi)(v). \]
For an automorphic character $\chi$ of $\U(W)$, its global theta lift $\Theta(\chi)$ is spanned by the automorphic forms
\[  \theta(\phi, \chi)(h) = \int_{[\U(W)]} \theta(\phi)(g,h) \cdot \overline{\chi(g)} \, dg \]
as $\phi$ ranges over elements of $\mathcal{S}( Y_{\A} \otimes W_{\A})$.
\vskip 5pt

(i) Recall the Borel subgroup $B = TU$ of $\U(V)$ which is the stabilizer of $X = E \cdot e$. For any character $\psi'$ of $U(k) \backslash U(\A)$, compute the $(U, \psi')$-Fourier coefficient
\[  \theta(\phi,\chi)_{U,\psi'}(h) = \int_{[U]} \theta(\phi,\chi)(uh) \cdot \overline{\psi'(u)}\, du. \]
\vskip 5pt

(ii) From your computation in (i), deduce that the global theta lift $\Theta(\chi)$ is nonzero for any $\chi$, and  is  cuspidal if and only if $\chi \ne \chi_V \circ i$ (see Problem 6 for the definition of $i$).
\vskip 5pt

(iii) (Challenging)  We would like to express the global theta lift $\theta(\phi, \chi_V \circ i)$ (which is noncuspidal) as an explicit Eisenstein series.
Observe that the map $\phi \mapsto \theta(\phi, \chi_V \circ i)$ is an equivariant map
\[  \Omega \twoheadrightarrow (\chi_V \circ i) \otimes  \Theta(\chi_V \circ i) \subset (\chi_V \circ i) \otimes \mathcal{A}(\U(V))  \]
and that
\[  \dim \Hom_{\U(W) \times \U(V)}(\Omega, (\chi_V \circ i) \otimes  \Theta(\chi_V \circ i) )  = 1. \]
Now we shall produce another element in this 1-dimensional vector space. 
\vskip 5pt

Recall from Problem 7(iv) that the map $\phi \mapsto (h \mapsto \Omega(h) \phi(0))$ defines an equivariant map
\[ j: \Omega \longrightarrow (\chi_V \circ i) \otimes  I(\chi_W) \]
whose image is isomorphic  to  $(\chi_V \circ i) \otimes  \Theta(\chi_V \circ i)$. 
Now the Eisenstein series is a $\U(V)$-equivariant map
\[  E(s,-) : I(s, \chi_W)  \longrightarrow \mathcal{A}(\U(V)) \]
defined by
\[  E(s, f)(h) = \sum_{\gamma \in B \backslash \U(V)}  f(\gamma g). \]
This converges only when $Re(s)$ is sufficiently large (actually $Re(s) >1/2$), but a basic theorem is that it admits a meromorphic continuation to $\C$ and that it is holomorphic at $s = 0$. Admitting this, we can thus consider $E (-) := E(0,-)$. 
\vskip 5pt

Now we have  the composite map 
\[ E \circ j: \Omega \longrightarrow (\chi_V \circ i) \otimes \mathcal{A}(\U(V)). \]
Show that this map is nonzero, so that its image is isomorphic to   $(\chi_V \circ i) \otimes  \Theta(\chi_V \circ i)$.
\vskip 5pt
Note however that the image of $E \circ j$ is not yet known to be equal to the submodule $\Theta(\chi_V \circ i) \subset \mathcal{A}(\U(V))$, but merely isomorphic to it. 
If we had known that the image is equal to $\Theta(\chi_V \circ i)$, then we would have immediately deduce that
 there is a nonzro constant $c \in \C^{\times}$ such that 
\[  \theta(\phi, \chi_V \circ i) = c \cdot E(j(\phi)). \]
Despite this, show  that the above identity holds (so that the image is indeed equal to $\Theta(\chi_V \circ i)$) (Hint: you may want to consider the Fourier expansion of both sides).
This is the so-called global Siegel-Weil formula.

\vskip 10pt

\noindent (11) The purpose of this exercise is to do global analog of  the very long Problem 8.  We shall use (the global analog of) the notations from Problem 8, as well as the same seesaw setup. 
Before that, let us explain the global seesaw identity in the context of a general seesaw in an ambient group $E$:
 \[
 \xymatrix{
  G_2  \ar@{-}[dr] \ar@{-}[d] & H_1  \ar@{-}[d] \\
  G_1  \ar@{-}[ur] & H_2}
\] 
We are now working over a number field $k$. Suppose $\Omega$ is the ``Weil representation" of $E$ and $\Omega$ is equipped with an automorphic realization
\[  \theta:  \Omega \longrightarrow \mathcal{A}(E). \]
For $\phi \in \Omega$ and $f \in \mathcal{A}(G_1)$, one has the global theta lift:
\[  \theta(\phi, f) (h) = \int_{[G_1]} \theta(\phi)(g,h) \cdot \overline{f(g)} \, dg \]
so that $\theta(\phi, f) \in \mathcal{A}(H_1)$. Likewise for $f' \in \mathcal{A}(H_2)$, one has $\theta(\phi, f') \in \mathcal{A}(G_2)$ defined by
\[   
\theta(\phi, f')  = \int_{[H_2]} \theta(\phi)(g,h) \cdot \overline{f'(h)} \, dh.\]
 Now the global seesaw identity is simply:
\[  \langle \theta(\phi, f), f' \rangle_{H_2} = \langle \theta(\phi, f'), f \rangle_{G_1}, \]
where we are using the Petersson inner products on $G_1$ and $H_2$ here.
Indeed, from definition, it follows that both sides are given by the double integral
\[  \int_{[G_1 \times H_2]} \theta(\phi)(g,h)  \cdot \overline{f(g)} \cdot  \overline{f'(h)} \, dg \, dh. \]
Hence the global seesaw identity is simply an application of Fubini's theroem: exchanging the order of integration. 
\vskip 5pt

Now we place ourselves in the context of Problem 8, with the follwing seesaw:
 \[
 \xymatrix{
  \U(V^{\square})  \ar@{-}[dr] \ar@{-}[d] & \U(W) \times \U(W)   \ar@{-}[d] \\
  \U(V) \times \U(V^-)  \ar@{-}[ur] & \U(W)^{\Delta}}
\]

 (i)  Taking the automorphic characters $f = \chi \otimes  {\chi'}^{-1} \chi_W|_{E^1}$ on $\U(V) \times \U(V^-)$ and $f' = \chi_V|_{E^1}$ on $\U(W^{\Delta})$, 
 write down the resulting global seesaw identity.
 
 \vskip 5pt
 
(ii) We now examine the RHS of the seesaw identity (the side of $W$'s). 
For $\phi \in \Omega_{V^-, W,\psi}$, show using Problem 8(iii) that
\[  \theta(\phi, {\chi'}^{-1}\chi_W|_{E^1})  \in \Theta_{V,W,\psi}(\chi') \cdot \chi_V|_{E^1}. \]
\vskip 5pt

(iii) For $\phi_1 \in \Omega_{V,W,\psi}$ and $\phi_2 \in \Omega_{V^-, W,\psi}$, show that
\[    \int_{[\U(W)]} \theta(\phi_1, \chi)(g) \cdot \theta(\phi_2,{\chi'}^{-1} \chi_W|_{E^1})(g) \cdot \overline{ \chi_V(g)}\, dg \] 
can be nonzero only if $\chi' = \chi$.
Moreover,  $\theta(\phi_1, \chi)$ is nonzero if and only if the above integral is nonzero for some $\phi_2$ (and taking $\chi' = \chi$).
\vskip 5pt

(iv)   The LHS of the global seesaw identity (the side of $V$'s) has the form:
\[   \int_{[\U(V) \times \U(V^-)]}  \theta( \phi_1,\otimes \phi_2, \chi_V|_{E^1}) \cdot  \overline{\chi(h_1)} \cdot \chi (h_2) \cdot \chi_W(h_2)^{-1} \, dh_1 \, dh_2, \]
where  
\[  \theta( \phi_1,\otimes \phi_2, \chi_V|_{E^1}) = \int_{[\U(W^{\Delta})]} \theta(\phi_1 \otimes \phi_2)(gh) \cdot  \overline{ \chi_V(g)}\, dg \]
is the global theta lift of $\chi_V|_{E^1}$ from $\U(W^{\Delta})$ to $\U(V^{\square})$.  As in Problem 8, we now need to explicate the theta lift 
$\theta( \phi_1,\otimes \phi_2, \chi_V|_{E^1})$. This is provided by the global Siegel-Weil formula of Problem 10(iii), which expresses $\theta( \phi_1,\otimes \phi_2, \chi_V|_{E^1})$ as an Eisenstein series $E(j(\phi_1 \otimes \phi_2))$ (where we recall that $j(\phi_1 \otimes \phi_2) \in I(\chi_W)$).

\vskip 5pt

(v)  On replacing $\theta( \phi_1,\otimes \phi_2, \chi_V|_{E^1})$ by the Eisenstein series $E(j(\phi_1 \otimes \phi_2))$, the integral in (iv) becomes (a special value of) the gloibal doubling zeta integral:
\[  Z(s, \chi)(f) = \int_{[\U(V) \times \U(V^-)]} E(s, f)(h_1, h_2) \cdot \overline{\chi(h_1)}^{-1} \cdot \chi(h_2) \cdot \overline{\chi_W(h_2)}^{-1} \, dh_1\, dh_2, \]
for $f_s \in I(s, \chi_W)$. The theory of this doubling zeta integral is discussed in Ellen Eischen's lectures. As discussed there, this doubling zeta integral represents the L-value $L(1/2 +s, \chi_E \cdot \chi_W^{-1})$.  More precisely, for $f_s \in I(s, \chi_W)$, we have:
\[  Z(s, \chi)(f_s) = \frac{L_E(\frac{1}{2} +s, \chi_E \cdot \chi_W^{-1})}{L(1+2s, \omega_{E/k})} \cdot \prod_v  \frac{Z_v(s,  \chi_v)(f_{s,v}) \cdot L(1+2s, \omega_{E_v/k_v})}{L_{E_v}(\frac{1}{2}+s, \chi_{E_v} \cdot \chi_{W,v}^{-1})}, \]
where the product over $v$ is finite (as almost all terns are equal to $1$). 
Using this identity and the last assertion in Problem 8(viii), prove Theorem 3.3 in the lecture notes.
\vskip 10pt

\noindent (12) Do the guided exercise in \S 3.13 of the lecture notes. This is
the global analog of Problem 9.

\vskip 15pt

\subsection{\bf Lecture 4}  $\text{}$
\vskip 10pt

\noindent (13) For an A-parameter $\psi$ considered in Lecture 4, i.e. 
Saito-Kurokawa type for $\PGSp_4$, Howe-PS type fop $\U_3$ and the short and long root type for $G_2$, compute the global component group $S_{\psi}$ and the quadratic character $\epsilon_{\psi}$. 
\vskip 15pt

\noindent (14) Do the same for the original Howe-PS A-parameter on $\PGSp_4$ defined as follows:
\[  \psi = \rho \times {\rm id} :  L_k \times \SL_2(\C) \longrightarrow \O_2(\C) \times \SL_2(\C)  \longrightarrow \Sp_4(\C), \]
where the second arrow is defined by the natural map associated to the tensor product of a 2-dimensional quadratic space and a 2-dimensional symplectic space (which yields 
a 4-dimensional symplectic space). In addition, note that a homomorphism 
\[  \rho: L_k \longrightarrow \O_2(\C) \]
determines:
\begin{itemize}
\item  by composition with the determinant map on $\O_2(\C)$ a quadratic \'etale $k$-algebra $E$
\item an automorphic character $\chi_{\rho}$ of $E^1 \cong E^{\times}/k^{\times}$.
\end{itemize}

 \vskip 15pt
 
 \noindent{\bf Acknowledgments:} I thank the organizers of the 2022 Arizona Winter School for their kind invitation to deliver a short course based on these lecture notes and for their hospitality and logistical support. I also thank the referees for their thorough work and pertinent questions. 
  The author is partially supported by a Singapore government MOE Tier 1 grant 
R-146-000-320-114.

\end{document}